\documentclass[a4paper, 11pt, final]{article}

\usepackage{amsfonts,enumerate}
\usepackage{enumerate}
\usepackage{graphicx,psfrag}
\newtheorem{proposition}{Proposition}[section]
\newtheorem{theorem}[proposition]{Theorem}
\newtheorem{lemma}[proposition]{Lemma}
\newtheorem{corollary}[proposition]{Corollary}
\newtheorem{definition}[proposition]{Definition}

\newenvironment{proof}{\smallskip\noindent\emph{\textbf{Proof.}}\hspace{1pt}}%
{\hspace{-5pt}{\nobreak\quad\nobreak\hfill\nobreak$\square$\vspace{8pt}%
\par}\smallskip\goodbreak}

\newenvironment{proofof}[1]{\smallskip\noindent\emph{\textbf{Proof of #1.}}%
\hspace{1pt}}{\hspace{-5pt}{\nobreak\quad\nobreak\hfill\nobreak%
$\square$\vspace{8pt}\par}\smallskip\goodbreak}

\newcommand{\Section}[1]{\section{#1}\setcounter{equation}{0}}

\newcommand{\pint}[1]{\mathaccent23{#1}}

\newcommand{\C}[1]{\mathbf{C^{#1}}}

\newcommand{\PC}{\mathbf{PC}}
\newcommand{\Cc}[1]{\mathbf{C_c^{#1}}}

\newcommand{\modulo}[1]{{\left|#1\right|}}
\newcommand{\norma}[1]{{\left\|#1\right\|}}
\newcommand{\Ref}[1]{{\rm(\ref{#1})}}
\newcommand{\reali}{{\mathbb{R}}}
\newcommand{\rsenza}{{\pint{\reali}^+}}

\newcommand{\naturali}{{\mathbb{N}}}
\newcommand{\tv}{\mathrm{TV}}
\newcommand{\BV}{\mathbf{BV}}

\renewcommand{\epsilon}{\varepsilon}
\renewcommand{\phi}{\varphi}
\renewcommand{\theta}{\vartheta}
\renewcommand{\L}[1]{\mathbf{L^#1}}
\newcommand{\W}[1]{\mathbf{W^{#1}}}

\title{Smooth and Discontinuous Junctions \\ in the $p$-System}

\author{Rinaldo M.~Colombo \\ \small Dipartimento di Matematica \\
  \small Universit\`a degli Studi di Brescia \\ \small Via Branze, 38
  \\ \small 25123 Brescia, Italy \\ \small \texttt{Rinaldo.Colombo@UniBs.it} \\
  \and Francesca Marcellini \\ \small Dip.~di Matematica e Applicazioni \\
  \small Universit\`a di Milano -- Bicocca \\ \small Via Cozzi, 53 \\
  \small 20126 Milano, Italy \\ \small
  \texttt{F.Marcellini@Campus.UniMiB.it}}

\begin{document}

\maketitle

\begin{abstract}

  \noindent Consider the $p$--system describing the subsonic flow of a
  fluid in a pipe with section $a = a(x)$. We prove that the resulting
  Cauchy problem generates a Lipschitz semigroup, provided the total
  variation of the initial datum and the oscillation of $a$ are
  small. An explicit estimate on the bound of the total variation of
  $a$ is provided, showing that at lower fluid speeds, higher total
  variations of $a$ are acceptable. An example shows that the bound on
  $\tv(a)$ is mandatory, for otherwise the total variation of the
  solution may grow arbitrarily.

  \medskip

  \noindent\textit{2000~Mathematics Subject Classification:} 35L65,
  76N10.

  \medskip


\end{abstract}

\Section{Introduction}
\label{sec:Intro}

Consider a gas pipe with smoothly varying section. In the isentropic
or isothermal approximation, the dynamics of the fluid in the pipe is
described by the following system of Euler equations:
\begin{equation}
  \label{eq:GL2}
  \left\{ 
    \begin{array}{l}
      \displaystyle
      \partial_{t} (a\rho) 
      + 
      \partial_{x} (a q) 
      = 
      0 
      \\ 
      \displaystyle
      \partial_{t} (aq) 
      +
      \partial_{x} 
      \left[ 
        a\left( \frac{q^{2}}{\rho } + p(\rho) \right)
      \right]  
      =
      p\left( \rho\right)  \, \partial_{x}a,
    \end{array}
  \right.  
\end{equation}
where, as usual, $\rho$ is the fluid density, $q$ is the linear
momentum density, $p = p(\rho)$ is the pressure and $a = a(x)$ is
\textit{cross-sectional area} of the tube. We provide a basic well
posedness result for (\ref{eq:GL2}), under the assumptions that the
initial data is subsonic, has sufficiently small total variation and
the oscillation in the pipe section $a=a(x)$ is also small. We provide
an explicit bound on the total variation of $a$. As it is physically
reasonable, as the fluid speed increases this bound decreases and
vanishes at sonic speed, see~(\ref{eq:M}).

As a tool in the study of~(\ref{eq:GL2}) we use the system recently
proposed for the case of a sharp discontinuous change in the pipe's
section between the values $a^-$ and $a^+$, see~\cite{BandaHertyKlar1,
  ColomboGaravelloSIAM, ColomboHertySachers}. This description is
based on the $p$-system
\begin{equation}
  \label{eq:HS}
  \!\left\{
    \!
    \begin{array}{l}
      \partial_{t} \rho + \partial_{x} q = 0
      \\ 
      \partial_{t} q + \partial_{x} 
      \left( \frac{q^{2}}{\rho } + p(\rho) \right) = 0
    \end{array}
  \right.
\end{equation}
equipped with a \textit{coupling condition} at the junction of the
form
\begin{equation}
  \label{eq:Junction}
  \Psi \left( a^-,(\rho,q)(t, 0-); a^+, (\rho,q)(t, 0+) \right) = 0
\end{equation}
whose role is essentially that of selecting stationary solutions.

Remark that the introduction of condition~(\ref{eq:Junction}) is
necessary as soon as the section of the pipe is not smooth.  The
literature offers different choices for this condition,
see~\cite{BandaHertyKlar1, ColomboGaravelloSIAM,
  ColomboHertySachers}. The construction below does not require any
specific choice of $\Psi$ in~(\ref{eq:Junction}), but applies to all
conditions satisfying minimal physically reasonable requirements,
see~\textbf{($\mathbf{\Sigma}$0)}--\textbf{($\mathbf{\Sigma}$2)}

On the contrary, if $a \in \W{1,1}$ the product in the right hand side
of the second equation in~(\ref{eq:GL2}) is well defined and
system~(\ref{eq:GL2}) is equivalent to the $2\times 2$ system of
conservation laws
\begin{equation}
  \label{eq:GL}
  \left\{ 
    \begin{array}{l}
      \partial_{t} \rho + \partial_{x} q = 
      -\frac{q}{a} \, \partial _{x} a 
      \\ 
      \partial_{t} q + \partial_{x} 
      \left( \frac{q^{2}}{\rho } + p(\rho) \right) =
      -\frac{q^{2}}{a\rho } \, \partial_{x}a \,.
    \end{array}
  \right.  
\end{equation}
Systems of this type were considered, for instance,
in~\cite{ChenHandbook, DalMasoLeflochMurat, GoatinLefloch, Hong,
  KroenerLefloch, LiuTransonic, Whitham}. In this case the stationary
solutions to~(\ref{eq:GL2}) are characterized as solutions to
\begin{equation}
  \label{eq:Stationary}
  \!\!\!
  \left\{
    \!\!\!\!
    \begin{array}{l}
      \partial_x (a(x) \, \hat q) = 0
      \\
      \partial_x \!
      \left( 
        \!
        a(x) \,\left( \frac{q^{2}}{\rho } + p(\rho) \right)
        \!
      \right) 
      =
      p(\hat\rho)\, \partial_x a
    \end{array}
  \right.
  \mbox{ or }
  \left\{ 
    \!\!\!\!
    \begin{array}{l}
      \partial_{x} q = 
      -\frac{q}{a} \, \partial _{x} a 
      \\ 
      \partial_{x} 
      \left( \frac{q^{2}}{\rho } + p(\rho) \right) =
      -\frac{q^{2}}{a\rho } \, \partial_{x}a \,,
    \end{array}
  \right.  
  \!\!\!
\end{equation}
see Lemma~\ref{lem:Weak} for a proof of the equivalence
between~(\ref{eq:GL}) and~(\ref{eq:GL2}).

Thus, the case of a smooth $a$ induces a unique choice for
condition~(\ref{eq:Junction}), see~(\ref{eq:Psi})
and~(\ref{eq:Sigma}). Even with this choice, in the case of the
isothermal pressure law $p(\rho) = c^2 \rho$, we show below that a
shock entering a pipe can have its strength arbitrarily magnified,
provided the total variation of the pipe's section is sufficiently
high and the fluid speed is sufficiently near to the sound speed, see
Section~\ref{sec:Multiple}. Recall, from the physical point of view,
that the present situation neglects friction, viscosity and the
conservation of energy. Moreover, this example shows the necessity of
a bound on the total variation of the pipe section in any well
posedness theorem for~(\ref{eq:GL2}).

The next section is divided into three parts, the former one deals
with a pipe with a single junction, the second with a pipe with a
piecewise constant section and the latter with a pipe having a
$\W{1,1}$ section. All proofs are gathered in Section~\ref{sec:Tech}.

\Section{Notation and Main Results}
\label{sec:Main}

Throughout this paper, $u$ denotes the pair $(\rho, q)$ so that, for
instance, $u^\pm = (\rho^\pm, q^\pm)$, $\bar u = (\bar\rho, \bar q)$,
$\ldots$. Correspondingly, we denote by $f(u) = \left( q,
  P(\rho,q)\right)$ the flow in~(\ref{eq:HS}).  Introduce also the
notation $\reali^+ = \left[0, +\infty \right[$, whereas $\rsenza =
\left]0, +\infty \right[$. Besides, we let $a(x\pm) = \lim_{\xi \to
  x\pm} a(\xi)$. Below, $B(u;\delta)$ denotes the open ball centered
in $u$ with radius~$\delta$.

The pressure law $p$ is assumed to satisfy the following requirement:
\begin{description}
\item[(P)] $p \in \C2(\reali^+; \reali^+)$ is such that for all $\rho
  > 0$, $p'(\rho) >0$ and $p''(\rho) \geq 0$.
\end{description}
\noindent The classical example is the $\gamma$-law, where $p(\rho) =
k \, \rho^\gamma$, for a suitable $\gamma \geq 1$.

Recall the expressions of the eigenvalues $\lambda_{1,2}$ and
eigenvectors $r_{1,2}$ of the $p$-system, with $c$ denoting the sound
speed,
\begin{equation}
  \label{eigenvectors}
  \begin{array}{c}
    \lambda_{1} (u)
    = 
    \frac{q}{\rho }-c\left( \rho\right) 
    \,,\qquad
    c\left( \rho\right) 
    =
    \sqrt{p^{\prime } (\rho)}
    \,,\qquad
    \lambda_{2} (u)
    = 
    \frac{q}{\rho }+c\left( \rho\right)
    \,,
    \\[5pt]
    r_{1} (u)
    =
    \left[ 
      \begin{array}{c}
        -1 
        \\ 
        -\lambda_{1} (u)
      \end{array}
    \right]
    \,,\qquad\qquad\qquad
    r_{2} (u)
    = 
    \left[ 
      \begin{array}{c}
        1 
        \\ 
        \lambda_{2} (u)
      \end{array}%
    \right]\,.
  \end{array}
\end{equation}
The \textit{subsonic} region is given by
\begin{equation}
  \label{subsonic2}
  A_{0} 
  =
  \left\{ 
    u \in \pint{\reali}^{+} \times \reali
    \colon
    \lambda_{1} (u) < 0 <\lambda_{2} (u) 
  \right\} \,.  
\end{equation}
For later use, we recall the quantities
\begin{eqnarray*}
  \mbox{flow of the linear momentum:}
  &
  \quad P(u) &=
  \frac{q^2}\rho+p(\rho)
  \\
  \nonumber
  \mbox{total energy density:}
  & 
  \quad E (u) &=
  \frac{q^2}{2\rho} + \rho \, \int_{\rho_*}^\rho \frac{p(r)}{r^2}\, dr
  \\
  \nonumber
  \mbox{flow of the total energy density:}
  & 
  \quad F (u) &=
  \frac{q}{\rho} \cdot \left( E(u) + p(\rho) \right) \,
\end{eqnarray*}
where $\rho_{*}>0$ is a suitable fixed constant. As it is well known,
see~\cite[formula~(3.3.21)]{DafermosBook}, the pair $(E,F)$ plays the
role of the (mathematical) entropy - entropy flux pair.

\subsection{A Pipe with a Single Junction}
\label{sec:Single}

This paragraph is devoted to~(\ref{eq:HS})--(\ref{eq:Junction}). Fix
the section $\bar a>\Delta$, with $\Delta>0$ and the state $\bar u \in
A_0$.

First, introduce a function $\Sigma = \Sigma(a^-,a^+;u^-)$ that
describes the effects of the junction when the section changes from
$a^-$ to $a^+$ and the state to the left of the junction is $u^-$. We
specify the choice of~(\ref{eq:Junction}) writing
\begin{equation}
  \label{eq:Psi}
  \Psi(a^-,u^-;a^+,u^+)
  =
  \left[  
    \begin{array}{c}
      a^{+} q^{+} - a^{-} q^{-} 
      \\ 
      a^{+} P (u^+) 
      -
      a^{-} P (u^-) 
    \end{array}
  \right]
  -
  \Sigma (a^-,a^+; u^- ) \,.
\end{equation}
We pose the following assumptions on $\Sigma$:
\begin{description}
\item[($\mathbf{\Sigma}$0)] $\Sigma \in \C1 \left( [\bar a-\Delta,
    \bar a+\Delta] \times B\left( \bar u;\delta \right); \reali^2
  \right)$.
\item[($\mathbf{\Sigma}$1)] $\Sigma (a,a;u^-) = 0$ for all $a \in
  [\bar a-\Delta, \bar a+\Delta]$ and all $u^- \in B ( \bar u;
  \delta)$.
\end{description}
\noindent Condition~\textbf{($\mathbf{\Sigma}$0)} is a natural
regularity condition. Condition~\textbf{($\mathbf{\Sigma}$1)} is aimed
to comprehend the standard \emph{``no junction''} situation: if $a^- =
a^+$, then the junction has no effects and $\Sigma$ vanishes.

Conditions~\textbf{($\mathbf{\Sigma}$0)}--\textbf{($\mathbf{\Sigma}$1)}
ensure the existence of stationary solutions to
problem~(\ref{eq:HS})--(\ref{eq:Junction}).

\begin{lemma}
  \label{lem:Stationary}
  Let~\textbf{($\mathbf{\Sigma}$0)}--\textbf{($\mathbf{\Sigma}$1)}
  hold. Then, for any $\bar a \in \rsenza$, $\bar u \in A_0$, there
  exists a positive $\bar\delta$ and a Lipschitz map
  \begin{equation}
    \label{eq:T}
    T \colon 
    \left]\bar a-\bar\delta, \bar a+ \bar \delta\right[
    \times 
    \left]\bar a-\bar\delta, \bar a+ \bar \delta\right[
    \times
    B \left( \bar u ; \bar\delta\right)
    \to
    A_0
  \end{equation}
  such that
  \begin{displaymath}
    \left.
      \begin{array}{l}
        \Psi (a^-,u^-; a^+,u^+) = 0
        \\
        a^- \in \left]\bar a-\bar\delta, \bar a+ \bar \delta\right[
        \\
        a^+ \in \left]\bar a-\bar\delta, \bar a+ \bar \delta\right[
        \\
        u^-, u^+ \in B \left( \bar u ; \bar\delta\right)
      \end{array}
    \right\}
    \iff
    u^+ = T (a^-,a^+;u^-) \,.
  \end{displaymath}
\end{lemma}

\noindent In particular, $T(\bar a, \bar a, \bar u) = \bar u$. We may
now state a final requirement on $\Sigma$:
\begin{description}
\item[($\mathbf{\Sigma}$2)] $\Sigma (a^-, a^0; u^-) + \Sigma \left(
    a^0,a^+;T(a^-,a^0; u^-) \right) = \Sigma (a^-,a^+; u^-)\,$.
\end{description}
With $T$ as in Lemma~\ref{lem:Stationary}. Alternatively,
by~(\ref{eq:Psi}), the above condition~\textbf{($\mathbf{\Sigma}$2)}
can be restated as
\begin{displaymath}
  \left.
    \begin{array}{rcl}
      \Psi(a^-,u^-; a^0,u^0) & = & 0
      \\
      \Psi(a^0,u^0; a^+,u^+) & = & 0
    \end{array}
  \right\}
  \Rightarrow
  \Psi(a^-,u^-;a^+,u^+) = 0 \,.
\end{displaymath}
Condition~\textbf{($\mathbf{\Sigma}$2)} says that if the two Riemann
problems with initial states $(a^-,u^-),(a^0,u^0)$ and
$(a^0,u^0),(a^+,u^+)$ both yield the stationary solution, then also
the Riemann problem with initial state $(a^-,u^-)$ and $(a^+,u^+)$ is
solved by the stationary solution.

Remark that the \emph{``natural''} choice~(\ref{eq:Sigma}) implied by
a smooth section satisfies~\textbf{($\mathbf{\Sigma}$0)},
\textbf{($\mathbf{\Sigma}$1)} and~\textbf{($\mathbf{\Sigma}$2)}.

Denote now by $\hat u$ a map satisfying
\begin{equation}
  \label{eq:Hat}
  \!\!\!\!\!
  \hat u (x)
  =
  \left\{
    \!\!\!
    \begin{array}{l@{\mbox{ if }}l}
      \hat u^-& x < 0
      \\
      \hat u^+ & x > 0
    \end{array}
  \right.
  \mbox{ with } 
  \begin{array}{l}
    \Psi \left( a^-, \hat u^-; a^+, \hat u^+ \right)
    =0,\!
    \\
    \hat u^-, \hat u^+ \in A_0 .
  \end{array}
\end{equation}
The existence of such a map follows from Lemma~\ref{lem:Stationary}.
Recall first the definition of weak $\Psi$-solution,
see~\cite[Definition~2.1]{ColomboGaravelloSIAM}
and~\cite[Definition~2.1]{ColomboHertySachers}.

\begin{definition}
  \label{def:WeakPsi}
  Let $\Sigma$
  satisfy~\textbf{($\mathbf{\Sigma}$0)}--\textbf{($\mathbf{\Sigma}$2)}. A
  weak $\Psi$-solution to~(\ref{eq:HS})--(\ref{eq:Junction}) is a map
  \begin{equation}
    \label{eq:regularity}
    \begin{array}{rcl}
      u 
      & \in & 
      \C0 \left( 
        \reali^+; 
        \hat u
        + 
        \L1 (\reali^{+}; \rsenza \times \reali) \right)
      \\ 
      u (t) & \in & \BV (\reali; \rsenza \times \reali)
      \quad \mbox{ for a.e. } t \in \reali^+
    \end{array}
  \end{equation}
  such that
  \begin{description}
  \item[(W)] for all $\phi \in \Cc1 (\rsenza \times \reali; \reali)$
    whose support does not intersect $x=0$
    \begin{displaymath}
      \int_{\reali^+} \int_{\reali} 
      \left( 
        u \, \partial_t \phi + 
        f(u) \, \partial_x \phi
      \right)
      \, dx \, dt = 0 \,;
    \end{displaymath}
  \item[($\mathbf{\Psi}$)] for a.e.~$t \in \reali^+$ and with $\Psi$
    as in~(\ref{eq:Psi}), the junction condition is met:
    \begin{displaymath}
      \Psi \left( a^-, u (t, 0-); a^+, u (t, 0+) \right) = 0 \,.
    \end{displaymath}
  \end{description}
  \noindent It is also an entropy solution if
  \begin{description}
  \item[(E)] for all $\phi \in \Cc1 (\rsenza \times \reali; \reali^+)$
    whose support does not intersect $x=0$
    \begin{displaymath}
      \int_{\reali^+} \! \int_{\reali} \left( E(u)
        \, \partial_t \phi + F(u) \, \partial_x \phi \right)
      \, dx \, dt \geq 0 \,.
    \end{displaymath}
  \end{description}
\end{definition}

\noindent In the particular case of a Riemann Problem,
i.e.~of~(\ref{eq:GL2}) with initial datum
\begin{displaymath}
  u (0,x) = \left\{
    \begin{array}{l@{\quad\mbox{ if }\quad}rcl}
      u^- & x & > & 0
      \\
      u^+ & x & < & 0 \,,
    \end{array}
  \right.
\end{displaymath}
Definition~\ref{def:WeakPsi} reduces
to~\cite[Definition~2.1]{ColomboHertySachers}.

To state the uniqueness property in the theorems below, we need to
introduce the following integral conditions,
following~\cite[Theorem~9.2]{BressanLectureNotes}, see
also~\cite[Theorem~8]{GuerraMarcelliniSchleper}
and~\cite{AmadoriGosseGuerra}. Given a function $u=u(t,x)$ and a point
$(\tau,\xi)$, we denote by $U^\sharp_{(u;\tau,\xi)}$ the solution of
the homogeneous Riemann Problem consisting
of~(\ref{eq:HS})--(\ref{eq:Junction})--(\ref{eq:Psi}) with initial
datum at time $\tau$
\begin{equation}
  \label{eq:sharp}
  w (\tau,x) = \left\{
    \begin{array}{l@{\quad\mbox{ if }\quad}rcl}
      \lim_{x\to \xi -} u(\tau,x) & x & < & \xi
      \\
      \lim_{x\to \xi +} u(\tau,x) & x & > & \xi \,.
    \end{array}
  \right.
\end{equation}
and with $\Sigma$ satisfying~\textbf{($\mathbf{\Sigma0}$)},
\textbf{($\mathbf{\Sigma1}$)} and~\textbf{($\mathbf{\Sigma2}$)}.
Moreover, define $U^\flat_{(u;\tau,\xi)}$ as the solution of the
linear hyperbolic Cauchy problem with constant coefficients
\begin{equation}
  \label{eq:flat}
  \left\{
    \begin{array}{l}
      \partial_t \omega + \partial_x \widetilde A\omega = 0
      \qquad t \geq \tau
      \\
      w(\tau,x) = u(\tau,x) \,,
    \end{array}
  \right.
\end{equation}
with $\widetilde A = D f \left( u(\tau,\xi) \right)$.

The next theorem applies~\cite[Theorem~3.2]{ColomboHertySachers}
to~\Ref{eq:HS} with the choice~(\ref{eq:Psi}) to construct the
semigroup generated
by~(\ref{eq:HS})--(\ref{eq:Junction})--(\ref{eq:Psi}). The uniqueness
part follows from~\cite[Theorem~2]{GuerraMarcelliniSchleper}.

\begin{theorem}
  \label{thm:Single}
  Let $p$ satisfy~\textbf{(P)} and $\Sigma$
  satisfy~\textbf{($\mathbf{\Sigma}$0)}--\textbf{($\mathbf{\Sigma}$2)}.
  Choose any $\bar a >0$, $\bar u \in A_0$. Then, there exist a
  positive $\Delta$ such that for all $a^-, a^+$ with
  $\modulo{a^--\bar a} < \Delta$ and $\modulo{a^+-\bar a} < \Delta$,
  there exist a map $\hat u$ as in~(\ref{eq:Hat}), positive $\delta,
  L$ and a semigroup $S \colon \reali^ + \times \mathcal{D} \to
  \mathcal{D}$ such that
  \begin{enumerate}
  \item $\mathcal{D} \supseteq \left\{ u \in \hat u + \L1(\reali; A_0)
      \colon \tv(u-\hat u) < \delta \right\}$.
  \item For all $u \in \mathcal{D}$, $S_0 u = u$ and for all $t,s \geq
    0$, $S_t S_s u = S_{s+t} u$.
  \item For all $u, u' \in \mathcal{D}$ and for all $t,t' \geq 0$,
    \begin{displaymath}
      \norma{S_t u - S_{t'} u'}_{\L1}
      \leq
      L \cdot \left(
        \norma{u - u'}_{\L1} + \modulo{t-t'}
      \right) \,.
    \end{displaymath}
  \item If $u \in \mathcal{D}$ is piecewise constant, then for $t$
    small, $S_t u$ is the gluing of solutions to Riemann problems at
    the points of jump in $u$ and at the junction at $x = 0$.
  \item For all $u \in \mathcal{D}$, the orbit $t \to S_t u$ is a weak
    $\Psi$-solution to~(\ref{eq:HS}).
  \item Let $\hat\lambda$ be an upper bound for the moduli of the
    characteristic speeds in $\bar B \left(\hat u(\reali),
      \delta\right)$. For all $u \in \mathcal{D}$, the orbit $u(t) =
    S_t u$ satisfies the integral conditions
    \begin{enumerate}[(\it i)]
    \item For all $\tau >0$ and $\xi \in \reali$,
      \begin{equation}
        \label{eq:first_int_cond}
        \lim_{h\rightarrow 0}
        \frac{1}{h}
        \int_{\xi-h\hat\lambda}^{\xi+h\hat\lambda}
        \norma{u (\tau+h, x) - U^\sharp_{(u;\tau,\xi)} (\tau + h,x)}
        \, dx
        = 0 \,.
      \end{equation}
    \item There exists a $C > 0$ such that for all $\tau > 0$, $a,b
      \in \reali$ and $\xi \in \left]a, b \right[$,
      \begin{equation}
        \label{eq:second_int_cond}
        \begin{array}{c}
          \displaystyle
          \frac{1}{h} 
          \int_{a+h\hat\lambda}^{b-h\hat\lambda}
          \norma{u (\tau+h, x) - U^\flat_{(u;\tau,\xi)}
            \left(\tau+h,x\right)} dx 
          \\
          \displaystyle
          \leq
          C
          \left[ \tv\left\{u(\tau);\left] a,b\right[ \right\} \right]^2.
        \end{array}
      \end{equation}
    \end{enumerate}
  \item If a Lipschitz map $w \colon \reali \to \mathcal{D}$
    satisfies~(\ref{eq:first_int_cond})--(\ref{eq:second_int_cond}),
    then it coincides with the semigroup orbit: $w(t) = S_t
    \left(w(0)\right)$.
  \end{enumerate}
\end{theorem}

\noindent The proof is deferred to Paragraph~\ref{sub:31}. Note that,
similarly to what happens in the standard case
of~\cite[Theorem~9.2]{BressanLectureNotes},
condition~(\ref{eq:second_int_cond}) is always satisfied at a
junction.

\subsection{A Pipe with Piecewise Constant Section}
\label{sec:Multiple}

We consider now a tube with piecewise constant section
\begin{displaymath}
  a 
  = 
  a_0 \, \chi_{]-\infty, x_1] }
  + 
  \sum_{j=1}^{n-1} a_j \, \chi_{[x_{j}, x_{j+1}[}
  +
  a_{n} \, \chi_{[x_n, +\infty[}
\end{displaymath}
for a suitable $n \in \naturali$. The fluid in each pipe is modeled
by~\Ref{eq:HS}. At each junction $x_j$, we require
condition~\Ref{eq:Junction}, namely
\begin{equation}
  \label{eq:Psii}
  \Psi (a_{j-1}, u_j^-; a_{j}, u_j^+) 
  =
  0
  \quad 
  \begin{array}{l}
    \mbox{ for all } j = 1, \ldots, n \mbox{, where}
    \\
    \displaystyle  u_j^\pm = \lim_{x \to x_j\pm} u_j (x) \,.
  \end{array}
\end{equation}
We omit the formal definition of $\Psi$-solution
to~(\ref{eq:HS})--(\ref{eq:Junction}) in the present case, since it is
an obvious iteration of Definition~\ref{def:WeakPsi}.

\begin{theorem}
  \label{thm:n}
  Let $p$ satisfy~\textbf{(P)} and $\Sigma$
  satisfy~\textbf{($\mathbf{\Sigma}$0)}--\textbf{($\mathbf{\Sigma}$2)}. For
  any $\bar a>0$ and any $\bar u \in A_0$ there exist positive $M,
  \Delta,\delta, L, \mathcal{M}$ such that for any profile satisfying
  \begin{description}
  \item[(A0)] $a \in \PC \left(\reali;\left]\bar a- \Delta, \bar
        a+\Delta\right[\right)$ with $\tv(a) < M$,
  \end{description}
  \noindent there exists a piecewise constant stationary solution
  \begin{displaymath}
    \hat u
    = 
    \hat u_0 \chi_{\left]-\infty, x_1\right[}
    +
    \sum_{j=1}^{n-1} \hat u_j \chi_{\left]x_{j},  x_{j+1}\right[}
    +
    \hat u_n \chi_{\left]x_{n}, +\infty\right[}
  \end{displaymath}
  to~(\ref{eq:HS})--(\ref{eq:Psii}) satisfying
  \begin{eqnarray}
    \nonumber
    & &
    \hat u_j \in A_0 \mbox{ with }
    \modulo{\hat u_j - \bar u} < \delta
    \mbox{ for }
    j=0, \ldots n
    \\
    \nonumber
    & &
    \Psi
    \left( 
      a_{j-1}, \hat u_{j-1}; a_{j}, \hat u_{j}
    \right)
    =
    0
    \mbox{ for }
    j = 1, \ldots, n
    \\
    \label{eq:tvhat}
    & &
    \tv(\hat u) \leq \mathcal{M} \, \tv(a)
  \end{eqnarray}
  and a semigroup $S^a \colon \reali^ + \times \mathcal{D}^a \to
  \mathcal{D}^a$ such that
  \begin{enumerate}
  \item $\mathcal{D}^a \supseteq \left\{ u \in \hat u + \L1(\reali;
      A_0) \colon \tv(u - \hat u) < \delta \right\}$.
  \item $S^{a}_0$ is the identity and for all $t,s \geq 0$, $S^a_t
    S^a_s = S^a_{s+t}$.
  \item For all $u, u' \in \mathcal{D}^a$ and for all $t, t' \geq 0$,
    \begin{displaymath}
      \norma{S^a_t u - S^a_{t'} u'}_{\L1}
      \leq
      L \cdot \left(
        \norma{(u) - u'}_{\L1} + \modulo{t-t'}
      \right) .
    \end{displaymath}
  \item If $u \in \mathcal{D}^a$ is piecewise constant, then for $t$
    small, $S_t u$ is the gluing of solutions to Riemann problems at
    the points of jump in $u$ and at each junction $x_j$.
  \item For all $u \in \mathcal{D}^a$, the orbit $t \to S^a_t u$ is a
    weak $\Psi$-solution to~(\ref{eq:HS})--(\ref{eq:Psii}).
  \item The semigroup satisfies the integral
    conditions~(\ref{eq:first_int_cond})--(\ref{eq:second_int_cond})
    in~6.~of Theorem~\ref{thm:Single}.
  \item If a Lipschitz map $w \colon \reali \to \mathcal{D}$
    satisfies~(\ref{eq:first_int_cond})--(\ref{eq:second_int_cond}),
    then it coincides with the semigroup orbit: $w(t) = S_t
    \left(w(0)\right)$.
  \end{enumerate}
\end{theorem}

\noindent Remark that $\delta$ and $L$ depend on $a$ only through
$\bar a$ and $\tv(a)$. In particular, all the construction above is
independent from the number of points of jump in $a$.  For every $\bar
u$, we provide below an estimate of $M$ at the leading order in
$\delta$ and $\Delta$, see~(\ref{eq:CK}) and~(\ref{eq:K}). In the case
of $\Sigma$ as in~(\ref{eq:Sigma}) and with the isothermal pressure
law, which obviously satisfies~\textbf{(P)},
\begin{equation}
  \label{eq:pIso}
  p(\rho) = c^2 \rho \,,
\end{equation}
the bounds~(\ref{eq:CK}) and~(\ref{eq:K}) reduce to the simpler
estimate
\begin{equation}
  \label{eq:M}
  M =
  \left\{
    \begin{array}{c@{\qquad\mbox{if }}rcl}
      \displaystyle
      \frac{\bar a}{4e} 
      & \bar v/c & \in & \left]0,1 / \sqrt2 \right],
      \\[8pt]
      \displaystyle
      \frac{\bar a}{4e} \frac{1-(\bar v/c)^2}{(\bar v/c)^2}
      & \bar v/c & \in & \left]1 / \sqrt2, 1\right[,
    \end{array}
  \right.
\end{equation}
where $\bar v = \bar q/\bar\rho$. Note that, as it is physically
reasonable, $M$ is a weakly decreasing function of $\bar v$, so that
at lower fluid speeds, higher values for the total variation of the
pipe's section can be accepted.

\begin{figure}[htpb]
  \centering
  \begin{psfrags}
    \psfrag{x}{$x$} \psfrag{da}{$\Delta a$} \psfrag{a}{$a$}
    \psfrag{s21}{$\sigma_2^-$} \psfrag{s22}{$\sigma_2^+$}
    \psfrag{s23}{$\sigma_2^{++}$} \psfrag{t}{$t$} \psfrag{up}{$u^+$}
    \psfrag{u}{$(u)$} \psfrag{2l}{$2l$} \psfrag{l}{$l$}
    \includegraphics[width=8cm]{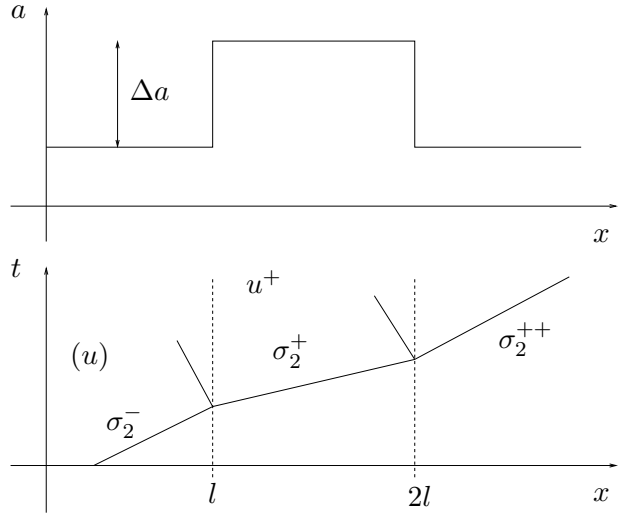}
  \end{psfrags}
  \caption{\label{fig:sugiu}A wave $\sigma_2^-$ hits a junction,
    giving rise to $\sigma_2^+$ which hits a second junction.}
\end{figure}

Furthermore, the estimates proved in Section~\ref{subs:22} show that
the total variation of the solution to~(\ref{eq:HS})--(\ref{eq:Psii})
may grow unboundedly if $\tv(a)$ is large. Consider the case in
Figure~\ref{fig:sugiu}.
A wave $\sigma_2^-$ hits a junction where the pipe's section increases
by $\Delta a>0$. From this interaction, the wave $\sigma_2^+$ of the
second family arises, which hits the second junction where the section
diminishes by $\Delta a$. At the leading term in $\Delta a$, we have
the estimate
\begin{eqnarray}
  \label{eq:sigma2++}
  \modulo{\sigma_2^{++}} 
  & \leq &
  \left(1 + \mathcal{K}(\bar v/c) \left(\frac{\Delta a}{a}\right)^2 \right) 
  \modulo{\sigma_2^-} \,,
  \qquad \mbox{ where}
  \\
  \label{eq:kgrande}
  \mathcal{K} (\xi)
  & = &
  \frac{-1 +8\,{\xi}^{2} -7\,{\xi}^{4} + 2\,{\xi}^{6}}{2(1-\xi)^3 \, (1+\xi)^3} \,,
\end{eqnarray}
see Section~\ref{subs:22} for the proof.  Note that $\mathcal{K}(0) =
-1$ whereas $\lim_{\xi \to 1-} \mathcal{K}(\xi) = +\infty$. Therefore,
for any fixed $\Delta a$, if $\bar v$ is sufficiently near to $c$,
repeating the interactions in Figure~\ref{fig:sugiu} a sufficient
number of times makes the $2$ shock waves arbitrarily large.

\subsection{A Pipe with a $\W{1,1}$ Section}
\label{sec:Smooth}

In this paragraph, the pipe's section $a$ is assumed to satisfy
\begin{displaymath}
  \mbox{\textbf{(A1)}}
  \hfill\hfill\qquad
  \left\{
    \begin{array}{l}
      a \in \W{1,1} \left(\reali;\left]\bar a-\Delta, \bar a+\Delta\right[
      \right) \mbox{ for suitable } \Delta >0, \ \bar a > \Delta
      \\
      \tv(a) < M \mbox{ for a suitable } M > 0
      \\
      a'(x) = 0 \mbox{ for a.e. } 
      x \in \reali \setminus [-X, X]
      \mbox{ for a suitable } X > 0 \,.
    \end{array}
  \right.
  \!\!
\end{displaymath}
For smooth solutions, the equivalence of~(\ref{eq:GL2})
and~(\ref{eq:GL}) is immediate. Note that the latter is in the
standard form of a 1D conservation law and the usual definition of
weak entropy solution applies, see for
instance~\cite[Definition~3.5.1]{SerreI}
or~\cite[Section~6]{CrastaPiccoli}. The definition below of weak
entropy solution to~(\ref{eq:GL2}) makes the two systems fully
equivalent also for non smooth solutions.

\begin{definition}
  \label{def:weak}
  A weak solution to~\Ref{eq:GL2} is a map
  \begin{displaymath}
    u \in \C0\left(
      \reali^+; \hat u+ \L1(\reali; \pint{\reali}^+ \times
      \reali) \right)
  \end{displaymath}
  such that for all $\phi \in \Cc1 (\pint{\reali}^+ \times \reali;
  \reali)$
  \begin{equation}
    \label{eq:weak}
    \int_{\reali^+} \! \int_{\reali}  \!
    \left( 
      \left[ \!
        \begin{array}{c}
          a\rho \\ a q
        \end{array}
        \! \right]
      \partial_t \phi + 
      \left[ \!
        \begin{array}{c}
          a q
          \\
          a P(u)
        \end{array}
        \! \right]
      \partial_x \phi
      +
      \left[ \!
        \begin{array}{c}
          0 \\ p(\rho) \partial_x a\!
        \end{array}
        \! \right]
      \, \phi
    \right)\!
    {dx} \, {dt}
    =
    0 \,.
  \end{equation}
  $u$ is an entropy weak solution if, for any $\phi \in \Cc1
  (\pint{\reali}^+ \times \reali; \reali)$, $\phi\geq0$,
  \begin{equation}
    \label{eq:Entropy}
    \int_{\reali^+} \! \int_{\reali}
    \left(
      a \, E(u) \, \partial_t \phi 
      +
      a \, F(u) \, \partial_x \phi 
    \right)
    \, dx \, dt
    \geq 
    0 \,.
  \end{equation}
\end{definition}

\begin{lemma}
  \label{lem:Weak}
  Let $a$ satisfy~\textbf{(A1)}. Then, $u$ is a weak entropy solution
  to~(\ref{eq:GL2}) in the sense of Definition~\ref{def:weak}, if and
  only if it is a weak entropy solution of~(\ref{eq:GL}).
\end{lemma}

\noindent The proof is deferred to Section~\ref{subs:23}.

Now, the section $a$ of the pipe is sufficiently regular to select
stationary solutions as solutions to either of the
systems~(\ref{eq:Stationary}), which are equivalent by
Lemma~\ref{lem:Weak}. Hence, the smoothness of $a$ also singles out a
specific choice of $\Sigma$,
see~\cite[formula~(14)]{GuerraMarcelliniSchleper}.

\begin{proposition}
  \label{prop:Sigma}
  Fix $a^-,a^+ \in \left]\bar a - \Delta, \bar a + \Delta\right[$ and
  $u^- \in A_0$. Choose a function $a$ strictly monotone, in $\C1$,
  that satisfies~\textbf{(A1)} with $a(-X-) = a^-$ and $a(X+) =
  a^+$. Call $\rho = R^a (x; u^- )$ the $\rho$-component of the
  corresponding solution to either of the Cauchy
  problems~(\ref{eq:Stationary}) with initial condition $u(-X) =
  u^-$. Then,
  \begin{enumerate}
  \item the function
    \begin{equation}
      \label{eq:Sigma}
      \Sigma (a^-, a^+, u^- )
      =
      \left[
        \begin{array}{c}
          0
          \\
          \displaystyle 
          \int_{-X}^{X} p \left( R^a(x; u^-) \right) a'(x) \, dx
        \end{array}
      \right]
    \end{equation}
    satisfies~\textbf{($\mathbf{\Sigma}$0)}--\textbf{($\mathbf{\Sigma}$2)};
  \item if $\tilde a$ is a strictly monotone function satisfying the
    same requirements above for $a$, the corresponding map $\tilde
    \Sigma$ coincides with $\Sigma$.
  \end{enumerate}
\end{proposition}

\noindent The basic well posedness theorem in the present $\W{1,1}$
case is stated similarly to Theorem~\ref{thm:n}.

\begin{theorem}
  \label{thm:W11}
  Let $p$ satisfy~\textbf{(P)}. For any $\bar a>0$ and any $\bar u \in
  A_0$ there exist positive $M, \Delta,\delta, L$ such that for any
  profile $a$ satisfying~\textbf{(A1)} there exists a stationary
  solution $\hat u$ to~(\ref{eq:GL2}) satisfying
  \begin{eqnarray*}
    & &
    \hat u \in A_0 \mbox{ with }
    \norma{\hat u(x) - \bar u} < \delta
    \mbox{ for all } x \in \reali
  \end{eqnarray*}
  and a semigroup $S^{a} \colon \reali^ + \times \mathcal{D}^{a} \to
  \mathcal{D}^{a}$ such that
  \begin{enumerate}
  \item $\mathcal{D}^{a} \supseteq \left\{ u \in \hat u + \L1(\reali;
      A_0) \colon \tv(u - \hat u) < \delta \right\}$.
  \item $S^{a}_0$ is the identity and for all $t,s \geq 0$, $S^{a}_t
    S^{a}_s = S^{a}_{s+t}$.
  \item for all $u, u' \in \mathcal{D}^a$ and for all $t,t' \geq 0$,
    \begin{displaymath}
      \norma{S^{a}_tu - S^{a}_{t'} u'}_{\L1}
      \leq
      L \cdot \left(
        \norma{u - u'}_{\L1} + \modulo{t-t'}
      \right) .
    \end{displaymath}
  \item for all $u \in \mathcal{D}^{a}$, the orbit $t \to S^{a}_t u$
    is solution to~(\ref{eq:HS}) in the sense of
    Definition~\ref{def:weak}.
  \item Let $\hat\lambda$ be an upper bound for the moduli of the
    characteristic speeds in $\bar B \left(\hat u(\reali),
      \delta\right)$. For all $u \in \mathcal{D}$, the orbit $u(t) =
    S_t u$ satisfies the integral conditions
    \begin{enumerate}[(\it i)]
    \item For all $\tau < 0$ and $\xi \in \reali$,
      \begin{equation}
        \label{eq:first_int_condBIS}
        \lim_{h\rightarrow 0}
        \frac{1}{h}
        \int_{\xi-h\hat\lambda}^{\xi+h\hat\lambda}
        \norma{u (\tau+h, x) - U^\sharp_{(u;\tau,\xi)} (\tau + h,x)}
        \, dx
        = 0 \,.
      \end{equation}
    \item There exists a $C > 0$ such that, for all $\tau > 0$, $a,b
      \in \reali$ and $\xi \in \left]a, b \right[$,
      \begin{equation}
        \label{eq:second_int_condBIS}
        \begin{array}{c}
          \displaystyle
          \frac{1}{h} 
          \int_{a+h\hat\lambda}^{b-h\hat\lambda}
          \norma{u (\tau+h, x) - U^\flat_{(u;\tau,\xi)}
            \left(\tau+h,x\right)} dx 
          \\[8pt]
          \displaystyle
          \leq
          C
          \left[ 
            \tv \left\{u(\tau);\left] a,b\right[ \right\} 
            +
            \tv \left\{ a; \left]a, b\right[ \right\}
          \right]^2.
        \end{array}
      \end{equation}
    \end{enumerate}
  \item If a Lipschitz map $w \colon \reali \to \mathcal{D}$
    solves~(\ref{eq:GL2}), then it coincides with the semigroup orbit:
    $w(t) = S_t \left(w(0)\right)$.
  \end{enumerate}
\end{theorem}

\noindent Thanks to Theorem~\ref{thm:n}, the proof is obtained
approximating $a$ with a piecewise constant function $a_n$. The
corresponding problems~(\ref{eq:HS})--(\ref{eq:Psii}) generate
semigroups defined on domains characterized by uniform bounds on the
total variation and with a uniformly bounded Lipschitz constants for
their time dependence. Then, we pass to the limit (see Section~4 for
the proof) and we follow the same procedure as
in~\cite[Theorem~9.2]{BressanLectureNotes} and~\cite[theorems~2
and~8]{GuerraMarcelliniSchleper} to characterize the solution.

As a byproduct of the proof of Theorem~\ref{thm:W11}, we also obtain
the following convergence result, relating the construction in
Theorem~\ref{thm:n} to that of Theorem~\ref{thm:W11}.

\begin{proposition}
  \label{prop:an}
  Under the same assumptions of Theorem~\ref{thm:W11}, for every $n
  \in \naturali$, choose a function $\beta_n$ such that:
  \begin{enumerate}[(i)]
  \item $\beta_n$ is piecewise constant with points of jump $y^1_n,
    \ldots, y^{m_n}_n$, with $y^1_n =-X$, $y^{m_n}_n = X$, and $\max_j
    ( y^{j+1}_n - y^j_n) \leq 1/n$.
  \item $\beta_n(x) = 0$ for all $x \in \reali \setminus [-X,X]$.
  \item $\beta_n \to a'$ in $\L1(\reali;\reali)$ with
    $\norma{\beta_n}_{\L1} \leq M$, with $M$ as in
    Theorem~\ref{thm:W11}.
  \end{enumerate}
  Define $\alpha_n(x) = a(-X-) + \int_{-X}^x \beta_n(\xi) \, d\xi$ and
  points $x^j_n \in \left]y^j_n, y^{j+1}_n\right[$ for $j= 1, \ldots,
  m_n-1$ and let
  \begin{displaymath}
    a_n
    = 
    a(-X-) \, \chi_{\strut]-\infty, x^1_j[}
    +
    \sum_{j=1}^{m_n-1} 
    \alpha_n(y^{j+1}_n) \, 
    \chi_{\strut [x^{j}_n, x^{j+1}_n[}
    +
    a(X+) \, \chi_{\strut [x^{m_n}_n, +\infty[}
  \end{displaymath}
  (see Figure~\ref{fig:alpha}) . Then, $a_n$ satisfies~\textbf{(A0)}
  and the corresponding semigroup $S^n$ constructed in
  Theorem~\ref{thm:n} converges pointwise to the semigroup $S$
  constructed in Theorem~\ref{thm:W11}.
\end{proposition}

\Section{Technical Proofs}
\label{sec:Tech}

\subsection{Proofs Related to Section~\ref{sec:Single}}
\label{sub:31}

The following equalities will be of use below:
\begin{equation}
  \label{eq:Useful}
  \partial_\rho P = - \lambda_1 \, \lambda_2
  \qquad\mbox{ and }\qquad
  \partial_q P = \lambda_1 + \lambda_2 \,.
\end{equation}

\begin{proofof}{Lemma~\ref{lem:Stationary}}
  Apply the Implicit Function Theorem to the equality $\Psi=0$ in a
  neighborhood of $(\bar a, \bar u, \bar a, \bar u)$, which satisfies
  $\Psi = 0$ by~\textbf{($\mathbf{\Sigma}$1)}. Observe that
  $\partial_{u} \Sigma (a, a; u^-) = 0$
  by~\textbf{($\mathbf{\Sigma}$1)}. Using~(\ref{eq:Useful}), compute
  \begin{eqnarray*}
    \det \partial_{u^+}\Psi 
    (\bar a, \bar u,\bar a, \bar u)
    & = &
    \det
    \left[
      \begin{array}{cc}
        - \partial_{\rho^+} \Sigma_1
        &
        \bar a
        \\
        \bar a\, \partial_{\rho^+}P
        &
        \bar a\, \partial_{q^+}P
      \end{array}
    \right]
    \\
    & = &
    \det
    \left[
      \begin{array}{cc}
        0
        &
        \bar a
        \\
        \bar a \, \partial_{\rho^+}P
        &
        \bar a \, \partial_{q^+}P
      \end{array}
    \right]
    \\
    & = &
    \bar a^2 \, \lambda_1(\bar u) \, \lambda_2(\bar u)
    \\
    & \neq &
    0 \,,
  \end{eqnarray*}
  completing the proof.
\end{proofof}

\begin{proofof}{Theorem~\ref{thm:Single}}
  Let $\Delta$ be defined as in
  Lemma~\ref{lem:Stationary}. Assumption~\textbf{(F)}
  in~\cite[Theorem~3.2]{ColomboHertySachers} follows
  from~\textbf{(P)}, thanks to~(\ref{eigenvectors}) and to the
  choices~(\ref{subsonic2})--(\ref{eq:Hat}). We now verify
  condition~\cite[formula~(2.2)]{ColomboHertySachers}. Recall that
  $D_{u^-} \Sigma (\bar a, \bar a; \bar u) = 0$
  by~\textbf{($\mathbf{\Sigma}$1)}. Hence, using~(\ref{eq:Useful}),
  \begin{eqnarray*}
    & &
    \det \left[ 
      D_{u^-}\Psi (\bar a, \bar u; \bar a, \bar u) \cdot r_{1} (\bar u)
      \quad
      D_{u^+}\Psi (\bar a, \bar u; \bar a, \bar u) \cdot r_{2} (\bar u)
    \right]  
    \\
    & = &
    \det
    \left[  
      \begin{array}{l@{\quad}l}
        \bar a \lambda_{1} (\bar u)
        +
        \partial_{\rho^-} \Sigma_1 (\bar a, \bar a; \bar u)
        + 
        \lambda_{1} \partial_{q^-} \Sigma_1 (\bar a, \bar a; \bar u)
        & 
        \bar a \lambda_{2} (\bar u)
        \\
        \bar a \left(\lambda_{1} (\bar u) \right)^{2}
        +
        \partial_{\rho^-} \Sigma_2 (\bar a, \bar a; \bar u)
        + 
        \lambda_{1}^- \partial_{q^-} \Sigma_2 (\bar a, \bar a; \bar u)
        & 
        \bar a \left(\lambda_{2} (\bar u) \right)^{2}
      \end{array}
    \right]
    \\
    & = &
    \det 
    \left[  
      \begin{array}{l@{\qquad}l}
        \bar a \, \lambda_{1} (\bar u)
        & 
        \bar a \, \lambda_{2} (\bar u)
        \\
        \bar a \, \left(\lambda_{1} (\bar u) \right)^{2} 
        & 
        \bar a \, \left(\lambda_{2} (\bar u) \right)^{2}
      \end{array}
    \right]
    \\
    & = &
    \bar a^2 \,  \lambda_1 (\bar u)\,
    \lambda_2 (\bar u) \left(\lambda_2(\bar u) - \lambda_1(\bar u) \right)
    \\
    & \neq &
    0 \,.
  \end{eqnarray*}
  The proof of 1.--5.~is completed
  applying~\cite[Theorem~3.2]{ColomboHertySachers}. The obtained
  semigroup coincides with that constructed
  in~\cite[Theorem~2]{GuerraMarcelliniSchleper}, where the uniqueness
  conditions~6.~and~7.~are proved.
\end{proofof}

\subsection{Proofs Related to Section~\ref{sec:Multiple}}
\label{subs:22}

We now work towards the proof of Theorem~\ref{thm:n}. We first use the
wave front tracking technique to construct approximate solutions to
the Cauchy problem~(\ref{eq:HS})--(\ref{eq:Psii}) adapting the wave
front tracking technique introduced
in~\cite[Chapter~7]{BressanLectureNotes}.

Fix an initial datum $u_o \in \hat u + \L1(\reali; A_0)$ and an
$\epsilon >0$. Approximate $u_o$ with a piecewise constant initial
datum $u_o^\epsilon$ having a finite number of discontinuities and so
that $\lim_{\epsilon \to 0} \norma{u_o^\epsilon - u_o}_{\L1} =
0$. Then, at each junction and at each point of jump in $u_o^\epsilon$
along the pipe, we solve the corresponding Riemann Problem according
to Definition~\ref{def:WeakPsi}. If the total variation of the initial
datum is sufficiently small, then Theorem~\ref{thm:Single} ensures the
existence and uniqueness of solutions to each Riemann Problem. We
approximate each rarefaction wave with a rarefaction fan, i.e.~by
means of (non entropic) shock waves traveling at the characteristic
speed of the state to the right of the shock and with size at most
$\epsilon$.

This construction can be extended up to the first time $\bar t_1$ at
which two waves interact in a pipe or a wave hits the junction.  At
time $\bar t_1$ the functions so constructed are piecewise constant
with a finite number of discontinuities. At any subsequent interaction
or collision with the junction, we repeat the previous construction
with the following provisions:
\begin{enumerate}
\item no more than $2$ waves interact at the same point or at the
  junction;
\item a rarefaction fan of the $i$-th family produced by the
  interaction between an $i$-th rarefaction and any other wave, is
  \emph{not} split any further;
\item when the product of the strengths of two interacting waves falls
  below a threshold $\check \epsilon$, then we let the waves cross
  each other, their size being unaltered, and introduce a \emph{non
    physical} wave with speed $\hat\lambda$, with $\hat\lambda >
  \sup_{(u)} \lambda_2(u)$; see~\cite[Chapter~7]{BressanLectureNotes}
  and the refinement~\cite{BaitiJenssen}.
\end{enumerate}
\noindent We complete the above algorithm stating how Riemann Problems
at the junctions are solved. We use the same rules as
in~\cite[\S~4.2]{ColomboGaravelloSIAM}
and~\cite[\S~5]{ColomboHertySachers}. In particular, at time $t=0$ and
whenever a physical wave with size greater than $\check \epsilon$ hits
the junction, the accurate solver is used, i.e.~the exact solution is
approximated replacing rarefaction waves with rarefaction fans. When a
non physical wave hits the junction, then we let it be refracted into
a non physical wave with the same speed $\hat\lambda$ and no other
wave is produced.

Repeating recursively this procedure, we construct a wave front
tracking sequence of approximate solutions $u_\epsilon$ in the sense
of~\cite[Definition~7.1]{BressanLectureNotes}.

At interactions of waves in a pipe, we have the following classical
result.

\begin{lemma}
  \label{lem:standard}
  Consider interactions in a pipe. Then, there exists a positive $K$
  with the properties:
  \begin{enumerate}
  \item An interaction between the wave $\sigma_1^-$ of the first
    family and $\sigma_2^-$ of the second family produces the waves
    $\sigma_1^+$ and $\sigma_2^+$ with
    \begin{equation}
      \modulo{\sigma_1^+ - \sigma_1^-}
      +
      \modulo{\sigma_2^+ - \sigma_2^-}
      \leq
      K \cdot \modulo{\sigma_1^- \sigma_2^-} \,.
    \end{equation}
  \item An interaction between $\sigma_i'$ and $\sigma_i''$ both of
    the same $i$-th family produces waves of total size $\sigma_1^+$
    and $\sigma_2^+$ with
    \begin{eqnarray*}
      \modulo{\sigma_1^+ - (\sigma_1''+ \sigma_1')}
      +
      \modulo{\sigma_2^+}
      \leq
      K \cdot \modulo{\sigma_1' \sigma_1''}
      \qquad \mbox{ if } i=1 \,,
      \\
      \modulo{\sigma_1^+}
      +
      \modulo{\sigma_2^+ - (\sigma_2''+ \sigma_2')}
      \leq
      K \cdot \modulo{\sigma_2' \sigma_2''}
      \qquad \mbox{ if } i=2 \,.
    \end{eqnarray*}
  \item An interaction between the physical waves $\sigma_1^-$ and
    $\sigma_2^-$ produces a non physical wave $\sigma_3^+$, then
    \begin{displaymath}
      \modulo{\sigma_3^+}\leq K \cdot \modulo{\sigma_1^-\sigma_2^-}.
    \end{displaymath}
  \item An interaction between a physical wave $\sigma$ and a non
    physical wave $\sigma_3^-$ produces a physical wave $\sigma$ and a
    non physical wave $\sigma_3^+$, then
    \begin{displaymath}
      \modulo{\sigma_3^+}-\modulo{\sigma_3^-}\leq K \cdot
      \modulo{\sigma\sigma_3^-}.
    \end{displaymath}
  \end{enumerate}
\end{lemma}

\noindent For a proof of this result
see~\cite[Chapter~7]{BressanLectureNotes}. Differently from the
constructions in~\cite{ColomboGaravelloSIAM, ColomboHertySachers}, we
now can not avoid the interaction of non physical waves with
junctions. Moreover, the estimates found therein do not allow to pass
to the limit $n \to +\infty$, $n$ being the number of junctions.

\begin{lemma}
  \label{lem:Junction}
  Consider interactions at the junction sited at $x_j$. There exist
  positive $K_1, K_2, K_3$ with the following properties.
  \begin{enumerate}
  \item \label{it:Inter1} The wave $\sigma_2^-$ hits the junction. The
    resulting waves $\sigma_1^+, \sigma_2^+$ satisfy\\
    \begin{minipage}{0.6\linewidth}
      \begin{eqnarray*}
        \modulo{\sigma_1^+}
        & \leq &
        K_1 \, \modulo{a_{j} - a_{j-1}} \modulo{\sigma_2^-} \,,
        \\
        \modulo{\sigma_2^+}
        &  \leq &
        \left( 1+ K_2 \, \modulo{a_{j} - a_{j-1}} \right)
        \modulo{\sigma_2^-}
        \\
        & \leq &
        e^{K_{2} \modulo{a_{j} - a_{j-1}}} \modulo{\sigma_2^-}\,.
      \end{eqnarray*}
    \end{minipage}
    \begin{minipage}{0.4\linewidth}
      \begin{center}
        \begin{psfrags}
          \psfrag{a}{$\sigma_2^-$} \psfrag{b}{$\sigma_1^+$}
          \psfrag{j}{$x_j$} \psfrag{c}{$\sigma_2^+$} \psfrag{u}{$\bar
            u$}
          \includegraphics[width=0.8\linewidth,
          height=3cm]{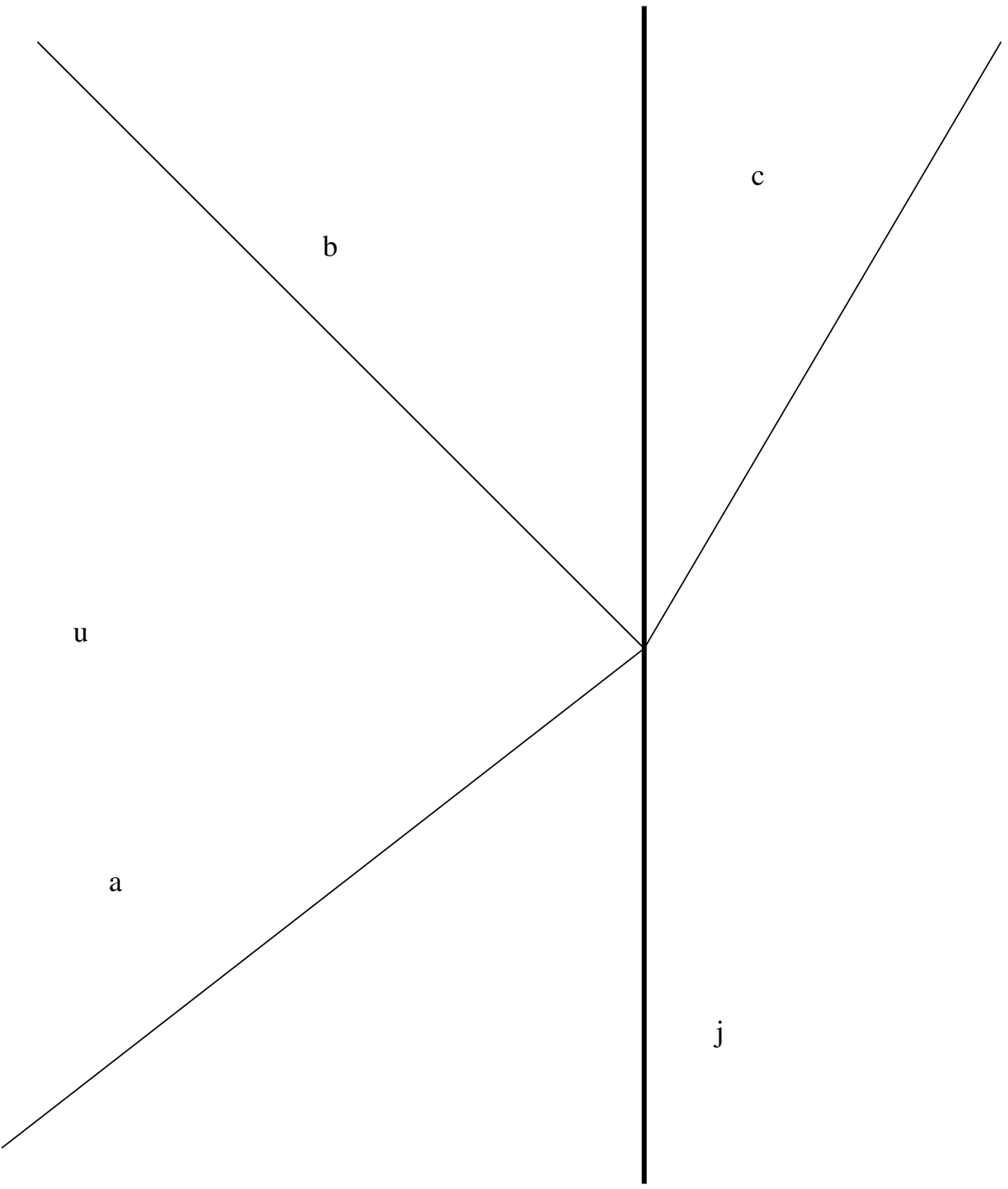}
        \end{psfrags}
      \end{center}
    \end{minipage}
  \item \label{it:Inter2} The non-physical wave $\sigma^-$ hits the
    junction. The resulting wave $\sigma^+$ satisfies\\
    \begin{minipage}{0.6\linewidth}
      \begin{eqnarray*}
        \modulo{\sigma^+}
        & \leq &
        \left( 1+ K_3 \, \modulo{a_{j} - a_{j-1}} \right)
        \modulo{\sigma^-}
        \\
        & \leq &
        e^{K_3 \modulo{a_{j} - a_{j-1}}} \modulo{\sigma^-}\,.
      \end{eqnarray*}
    \end{minipage}
    \begin{minipage}{0.4\linewidth}
      \begin{center}
        \begin{psfrags}
          \psfrag{a}{$\sigma^-$} \psfrag{j}{$x_j$}
          \psfrag{c}{$\sigma^+$} \psfrag{u}{$\bar u$}
          \includegraphics[width=0.8\linewidth,
          height=2.5cm]{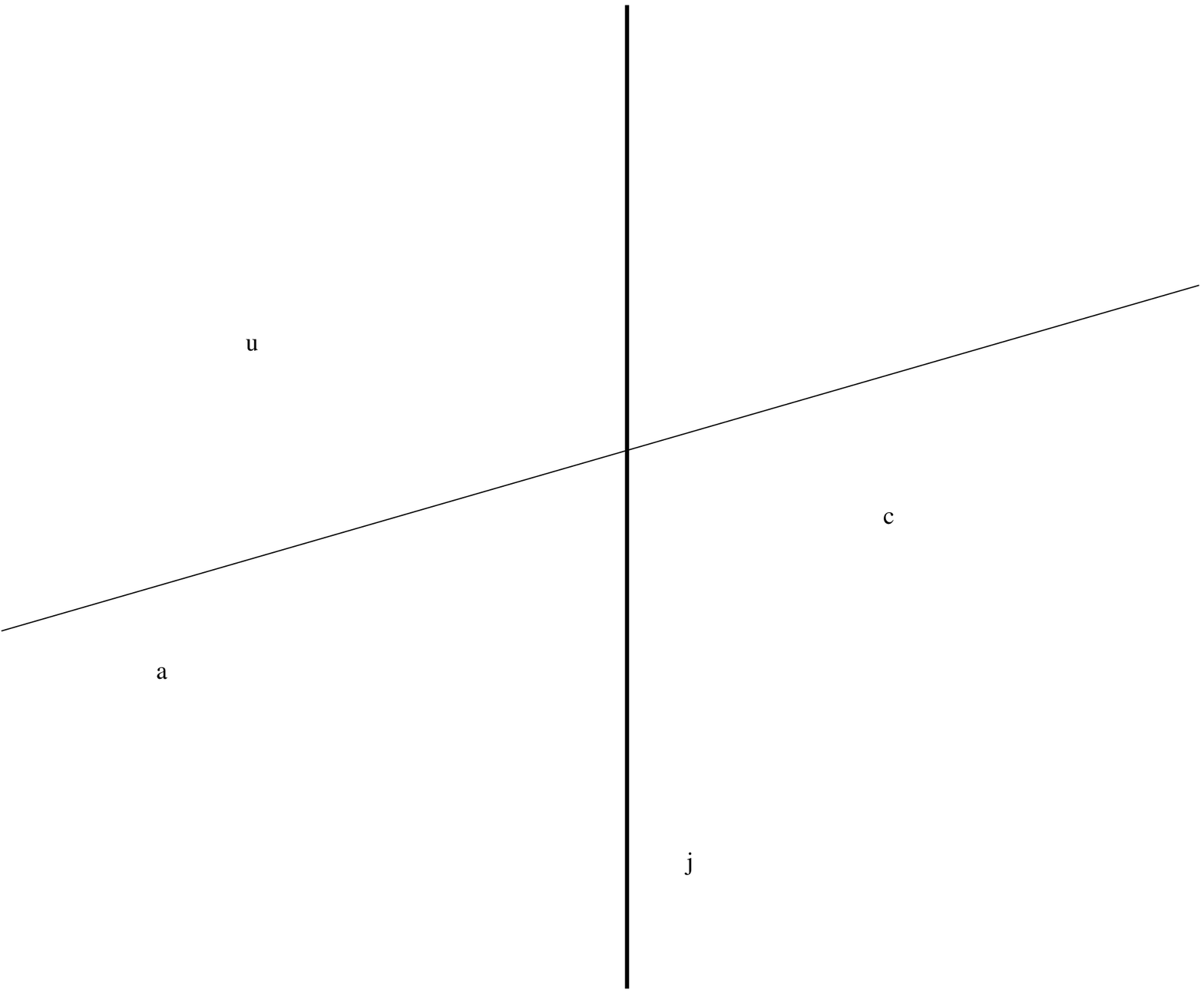}
        \end{psfrags}
      \end{center}
    \end{minipage}
  \end{enumerate}
\end{lemma}

\begin{proof}
  Use the notation in the figure above. Recall that $\sigma_1^+$ and
  $\sigma_2^+$ are computed through the Implicit Function Theorem
  applied to a suitable combination of the Lax curves
  of~(\ref{eq:HS}), see~\cite[Proposition~2.4]{ColomboGaravelloSIAM}
  and~\cite[Proposition~2.2]{ColomboHertySachers}. Repeating the proof
  of Theorem~\ref{thm:Single} one shows that the Implicit Function
  Theorem can be applied. Therefore, the regularity of the Lax curves
  and~\textbf{(P)} ensure that $\sigma_1^+ = \sigma_1^+
  \left(\sigma_2^-, a_{j} - a_{j-1}; \bar u\right)$ and $\sigma_2^+ =
  \sigma_2^+ \left(\sigma_2^-, a_{j} - a_{j-1}; \bar u\right)$. An
  application of~\cite[Lemma~2.5]{BressanLectureNotes}, yields
  \begin{eqnarray}
    \nonumber
    \!\!\!\!\!\!\!\!
    \left.\!\!\!
      \begin{array}{rcl}
        \sigma_1^+ \left(0, a_{j} - a_{j-1}; \bar u\right)
        & = &
        0
        \\
        \sigma_1^+ \left(\sigma_2^-, 0; \bar u\right)
        & = &
        0
      \end{array}
      \!\!\!
    \right\}
    \!\!\!& \Rightarrow&\!\!\!
    \modulo{\sigma_1^+}
    \leq
    K_1 \modulo{a_{j} - a_{j-1}} \, \modulo{\sigma_2^-} \,,
    \\
    \label{eq:type}
    \!\!\!\!\!\!\!\!
    \left.\!\!\!
      \begin{array}{rcl}
        \sigma_2^+ \left(0, a_{j} - a_{j-1}; \bar u\right)
        & = &
        0
        \\
        \sigma_2^+ \left(\sigma_2^-, 0; \bar u\right)
        & = &
        \sigma_2^-
      \end{array}
      \!\!\!
    \right\}
    \!\!\!& \Rightarrow &\!\!\!
    \modulo{\sigma_2^+ - \sigma_2^-}
    \leq
    K_2 \modulo{a_{j} - a_{j-1}} 
    \modulo{\sigma_2^-}
    \\
    \nonumber
    \!\!\!&\Rightarrow&\!\!\!
    \modulo{\sigma_2^+}
    \leq
    \left[ 1 + K_2 \modulo{a_{j} - a_{j-1}} \right] \! \modulo{\sigma_2^-} \,,
  \end{eqnarray}
  completing the proof of~1. The estimate at~2.~is proved similarly.
\end{proof}

We now aim at an improvement of~(\ref{eq:type}). Solving the Riemann
problem at the interaction in case~1.~amounts to solve the system
\begin{equation}
  \label{eq:PbR}
  \mathcal{L}_2 \left( 
    T \left( 
      \mathcal{L}_1 ( \bar u; \sigma_1^+) 
    \right)
    ; \sigma_2^+
  \right)
  = 
  T \left(
    \mathcal{L}_2(\bar u; \sigma_2^-)
  \right).
\end{equation}
By~(\ref{eigenvectors}), the first order expansions in the wave's
sizes of the Lax curves exiting $u$ are
\begin{displaymath}
  \mathcal{L}_1 (u; \sigma) 
  =
  \left[
    \begin{array}{c}
      \rho - \sigma + o(\sigma)
      \\
      q - \lambda_1(u) \sigma + o(\sigma)
    \end{array}
  \right]
  \mbox{ and }\ 
  \mathcal{L}_2 (u; \sigma) 
  = 
  \left[
    \begin{array}{c}
      \rho + \sigma + o(\sigma)
      \\
      q + \lambda_2(u) \sigma + o(\sigma)
    \end{array}
  \right],
\end{displaymath}
while the first order expansion in the size's difference $\Delta a =
a^+ - a^-$ of the map $T$ defined at~(\ref{eq:T}), with $v = q/\rho$,
is
\begin{eqnarray}
  \label{eq:TT}
  T(a, a+\Delta a; u)
  & = &
  \left[
    \begin{array}{c}
      \left( 
        1 + H
        \, \frac{\Delta a}{a}
      \right)
      \rho
      + o(\Delta a)
      \\
      \left( 1- \frac{\Delta a}{a} \right) q
      + o(\Delta a)
    \end{array}
  \right]
  \mbox{, \ where}
  \\
  \nonumber
  H 
  & = &
  \frac{v^2 + \frac{\partial_{a^+}\Sigma - p(\rho)}{\rho}}{c^2 - v^2} \,.
\end{eqnarray}
Inserting these expansions in~(\ref{eq:PbR}), we get the following
linear system for $\sigma_1^+, \sigma_2^+$:
\begin{displaymath}
  \left\{
    \begin{array}{rcrcr}
      \displaystyle
      -\left( 1 + \bar H \frac{\Delta a}{\bar a} \right) \sigma_1^+
      & + &
      \sigma_2^+
      & = &
      \left( 1 + \bar H \frac{\Delta a}{\bar a} \right) \sigma_2^-
      \\
      \displaystyle
      - \left( 1 - \frac{\Delta a}{\bar a} \right) \bar \lambda_1 
      \, \sigma_1^+
      & + & \left( 1 + \bar G \frac{\Delta a}{\bar a} \right) \bar \lambda_2
      \, \sigma_2^+
      & = &
      \left( 1 - \frac{\Delta a}{\bar a} \right) \bar \lambda_2 \sigma_2^-
    \end{array}
  \right.
\end{displaymath}
where
\begin{displaymath}
  \bar H 
  = 
  \frac{\bar v^2 + 
    \left( 
      \partial_{a^+}\Sigma (\bar a, \bar a, \bar u) -
      p(\bar \rho) \right)/ \bar \rho}{c^2-\bar v^2}
  \quad \mbox{ and } \quad
  \bar G 
  =
  \frac{(c'(\bar\rho) \bar\rho - \bar v)\bar H - \bar v}{\bar v +c}
\end{displaymath}
and all functions are computed in $\bar u$. The solution is
\begin{eqnarray}
  \label{eq:sigma1}
  \sigma_1^+
  & = &
  -\frac{\bar \lambda_2}{2c} \, ( 1 + \bar G + \bar H) \, \frac{\Delta a}{a}
  \, \sigma_2^-
  \\
  \label{eq:sigma2}
  \sigma_2^+
  & = &
  \left(
    1
    -
    \frac{\bar\lambda_1 \bar H + \bar \lambda_2 (1+\bar G)}{2c}
    \, \frac{\Delta a}{a}
  \right)
  \sigma_2^-
\end{eqnarray}
which implies the following first order estimate for the coefficients
in the interaction estimates of Lemma~\ref{lem:Junction}:
\begin{equation}
  \label{eq:K}
  \!\!\!
  \begin{array}{rcl}
    K_1
    & = &
    \displaystyle
    \frac{1}{2a}
    \modulo{
      \frac{
        1
        +
        \frac{c'\rho}{c} \, \left(\frac{v}{c}\right)^2 
        +
        \frac{1}{c^2} \left( \frac{c'\rho}{c}+1\right)
        \frac{\partial_{a^+}\Sigma - p(\rho)}{\rho}    
      }{1 - \left(\frac{v}{c}\right)^2}} \,,
    \\
    K_2
    & = &
    \displaystyle
    \frac{1}{2a}
    \modulo{
      \frac{
        1
        -
        2 \left( \frac{v}{c} \right)^2 
        +
        \frac{c'\rho}{c} \left( \frac{v}{c} \right)^2 
        +
        \frac{1}{c^2} \left( \frac{c'\rho}{c}-1\right)
        \frac{\partial_{a^+}\Sigma - p(\rho)}{\rho}    
      }{1 - \left(\frac{v}{c}\right)^2}} \,.
  \end{array}
\end{equation}

The estimate~(\ref{eq:sigma2}) directly implies the following corollary.

\begin{corollary}
  If $\modulo{a_j - a_{j-1}}$ is sufficiently small, then $\sigma_2^+$
  and $\sigma_2^-$ are either both rarefactions or both shocks.
\end{corollary}

Denote by $\sigma^j_{i,\alpha}$ the wave belonging to the $i$-th
family and sited at the point of jump $x^\alpha$, with $x^\alpha$ in
the $j$-th pipe $I_j$, where we set $I_0 = \left]-\infty, x_1\right[$,
$I_j =\left]x_j, x_{j+1}\right[$ for $j=1, \ldots, n-1$ and $I_n =
\left]x_n, +\infty \right[$. Aiming at a bound on the Total Variation
of the approximate solution, we define the Glimm-like functionals,
see~\cite[formul\ae~(7.53) and~(7.54)]{BressanLectureNotes} or
also~\cite{DafermosBook, Glimm, Lax1973, SmollerBook},
\begin{eqnarray}
  \nonumber
  V
  & = &
  \sum_{j=0}^{n}
  \sum_{x^\alpha \in I_j}
  \left(
    \modulo{\sigma^j_{1,\alpha}} e^{C\sum_{h=1}^j \modulo{a_h-a_{h-1}}} 
    +
    \modulo{\sigma^j_{2,\alpha}} e^{C\sum_{h=j}^{n-1} \modulo{a_{h+1}-a_h}} 
  \right)
  \\
  \nonumber
  & &
  +
  \sum_{j=0}^{n}
  e^{C\sum_{h=j}^{n-1} \modulo{a_{h+1}-a_h}}
  \sum_{\sigma\mbox{ non physical in }I_j} \modulo{\sigma} \,,
  \\
  \nonumber
  Q 
  & = &
  \sum_{(\sigma_{i,\alpha}^j, \sigma_{i',\alpha'}^{j'}) \in \mathcal{A}}
  \modulo{\sigma_{i,\alpha}^j \sigma_{i',\alpha'}^{j'}} \,,
  \\
  \label{eq:upsilon}
  \Upsilon
  & = &
  V + Q \,,
\end{eqnarray}
where $C$ is a positive constant to be specified below.  $\mathcal{A}$
is the set of pairs $(\sigma_{i,\alpha}^j, \sigma_{i',\alpha'}^{j'})$
of approaching waves, see~\cite[Paragraph~3,
Section~7.3]{BressanLectureNotes}. The $i$-wave $\sigma_{i,\alpha}^j$
sited at $x_\alpha$ and the $i'$-wave $\sigma_{i',\alpha'}^{j'}$ sited
at $x_{\alpha'}$ are approaching if either $i<i'$ and $x_\alpha >
x_{\alpha'}$, or if $i=i'<3$ and $\min\{\sigma^j_{i,\alpha},
\sigma^{j'}_{i',\alpha'}\} < 0$, independently from $j$ and $j'$.  As
usual, non physical waves are considered as belonging to a fictitious
linearly degenerate $3$rd family, hence they are approaching to all
physical waves to their right.

It is immediate to note that the weights $\exp \left( C\sum_{h=1}^j
  \modulo{a_h-a_{h-1}} \right)$ and $\exp \left( C\sum_{h=j}^{n-1}
  \modulo{a_{h+1}-a_h} \right)$ in the definition of $V$ are uniformly
bounded:
\begin{equation}
  \label{eq:BoundedWeights}
  \forall \, j
  \qquad
  \left\{
    \begin{array}{ccccc}
      1 
      & \leq & \exp \left(C\sum_{h=1}^j \modulo{a_h-a_{h-1}}\right)
      & \leq & \exp \left( C\, \tv (a) \right) ,
      \\
      1 
      & \leq & \exp \left( C\sum_{h=j}^{n-1} \modulo{a_{h+1}-a_h}\right)
      & \leq & \exp \left( C\, \tv (a) \right).
    \end{array}
  \right.
\end{equation}

Below, the following elementary inequality is of use: if $a < b$, then
$e^a - e^b < -(b-a) e^a$.
\begin{lemma}
  \label{lem:Upsilon}
  There exists a positive $\delta$ such that if an
  $\epsilon$-approximate wave front tracking solution $u = u(t,x)$ has
  been defined up to time $\bar t$, $\Upsilon \left( u (\bar t-)
  \right) < \delta$ and an interaction takes place at time $\bar t$,
  then the $\epsilon$-solution can be extended beyond time $\bar t$
  and $\Upsilon \left( u (\bar t+) \right) < \Upsilon \left( u (\bar
    t-) \right)$.
\end{lemma}

\begin{proof}
  Thanks to~(\ref{eq:BoundedWeights}) and Lemma~\ref{lem:standard},
  the standard interaction estimates,
  see~\cite[Lemma~7.2]{BressanLectureNotes}, ensure that $\Upsilon$
  decreases at any interaction taking place in the interior of $I_j$,
  for any $j= 0, \ldots, n$.

  Consider now an interaction at $x_j$. In the case of~\ref{it:Inter1}
  in Lemma~\ref{lem:Junction},
  \begin{eqnarray*}
    & &
    \Delta Q
    \\
    & \leq &
    \sum_{(\sigma_1^+, \sigma_{i,\alpha}) \in \mathcal{A}}
    \modulo{\sigma_1^+ \sigma_{i,\alpha}}
    +
    \sum_{(\sigma_2^+, \sigma_{i,\alpha}) \in \mathcal{A}}
    \modulo{\sigma_{i,\alpha}} 
    \left( \modulo{\sigma_2^+ } - \modulo{\sigma_2^- } \right)
    \\
    & \leq &
    \left(
      K_1 \, \modulo{a_j - a_{j-1}} 
      \sum_{i,\alpha} \modulo{\sigma_{i,\alpha}}
      +
      \left( e^{K_2 \modulo{a_j - a_{j-1}}} - 1 \right)
      \sum_{i,\alpha} \modulo{\sigma_{i,\alpha}}
    \right)
    \modulo{\sigma_2^-}
    \\
    & \leq &
    (K_1+K_2)\, \Upsilon(\bar t -) 
    \, \modulo{a_j - a_{j-1}} \, \modulo{\sigma_2^-}
    \\
    & \leq &
    (K_1+K_2) \, \delta \, \modulo{a_j - a_{j-1}} \, \modulo{\sigma_2^-} \,.
    \\
    & &
    \Delta V
    \\
    & \leq &
    e^{C\sum_{h=1}^{j-1} \modulo{a_h - a_{h-1}}} \modulo{\sigma_1^+} 
    + 
    e^{C\sum_{h=j}^{n-1} \modulo{a_{h+1}-a_h}}\modulo{\sigma_2^+}
    -
    e^{C\sum_{h=j-1}^{n-1} \modulo{a_{h+1}-a_h}}\modulo{\sigma_2^-}
    \\
    & \leq &
    e^{C\sum_{h=1}^{j-1} \modulo{a_h-a_{h-1}}} 
    \left(  K_1 \modulo{ a_j - a_{j-1}} \modulo{\sigma_2^-}\right) 
    \\
    & + &
    \left( 
      e^{C\sum_{h=j}^{n-1} \modulo{a_{h+1}-a_h}} 
      e^{K_2 \modulo{ a_j - a_{j-1}}} 
      -
      e^{C\sum_{h=j-1}^{n-1} \modulo{a_{h+1}-a_h}}
    \right)
    \modulo{\sigma_2^-}
    \\
    &\leq&
    \left(
      K_1 \vert a_j - a_{j-1}\vert
      e^{C\sum_{h=1}^{j-1} \modulo{a_h-a_{h-1}}}
    \right) \, \modulo{\sigma_2^-}
    \\
    & +& 
    e^{C\sum_{h=j}^{n-1} \modulo{a_{h+1}-a_h}}
    \left( 
      e^{K_2\vert a_j - a_{j-1}\vert}  
      - 
      e^{C\modulo{a_j - a_{j-1}} }  
    \right)
    \, \modulo{\sigma_2^-}
    \\
    & \leq &
    \left(
      K_1 \modulo{ a_j - a_{j-1}}
      e^{C\sum_{h=1}^{j-1} \modulo{a_h-a_{h-1}}}
    \right) \, \modulo{\sigma_2^-}
    \\
    & & 
    -
    (C-K_2) \modulo{ a_j - a_{j-1}} \, e^{K_2\modulo{a_j - a_{j-1}}}
    \, e^{C\sum_{h=j}^{n-1} \modulo{a_{h+1}-a_h}}
    \, \modulo{\sigma_2^-}
    \\
    & \leq &
    \left( 
      ( K_1+K_2) \left( 1 + e^{K_2\modulo{a^+- a^-}} \right) 
      e^{C\tv(a)} - C
    \right) \,
    \modulo{a_j - a_{j-1}} \, \modulo{\sigma_2^-}.
    \\
    & &
    \Delta \Upsilon
    \\
    & \leq &
    \left( 
      (K_1+K_2) 
      \left( 1 + e^{K_2\modulo{a^+- a^-}} \right) e^{C\tv(a)} 
      + 
      (K_1+K_2) \delta - C 
    \right)
    \modulo{a_j - a_{j-1}} \, \modulo{\sigma_2^-}.
  \end{eqnarray*}
  Choosing now, for instance,
  \begin{equation}
    \label{eq:CK}
    \!\!\!\!\!
    \delta < 1
    \,, \;\;
    C = \frac{1}{\tv(a)}
    \,,\;\;
    \modulo{a^+ - a^-} \leq \frac{\ln 2}{K_{2}}
    \;\mbox{ and }\;
    \tv(a) < \frac{1}{4(K_1+K_2)e} 
  \end{equation}
  the monotonicity of $\Upsilon$ in this first case is proved.

  Consider an interaction as in~\ref{it:Inter2}. of
  Lemma~\ref{lem:Junction}. Then, similarly,
  \begin{eqnarray*}
    \Delta Q
    & \leq &
    \sum_{(\sigma^+, \sigma_{i,\alpha}) \in \mathcal{A}}
    \modulo{\sigma_{i,\alpha}} 
    \left( \modulo{\sigma^+ } - \modulo{\sigma^- } \right)
    \\
    & \leq &
    \left( e^{K_{3} \modulo{a_j - a_{j-1}}} - 1 \right)
    \sum_{i,\alpha} \modulo{\sigma_{i,\alpha}}
    \,
    \modulo{\sigma^-}
    \\
    & \leq &
    K_{3}\, \Upsilon(\bar t -) 
    \, \modulo{a_j - a_{j-1}} \, \modulo{\sigma^-}
    \\
    & \leq &
    K_{3}\,\delta \, \modulo{a_j - a_{j-1}} \, \modulo{\sigma^-} \,.
    \\
    \Delta V
    & \leq &
    e^{C\sum_{h=j}^{n-1} \modulo{a_{h+1}-a_h}}\modulo{\sigma^+}
    -
    e^{C\sum_{h=j-1}^{n-1} \modulo{a_{h+1}-a_h}}\modulo{\sigma^-}
    \\
    & \leq &
    \left( 
      e^{C\sum_{h=j}^{n-1} \modulo{a_{h+1}-a_h}} 
      e^{K_{3}\modulo{ a_j - a_{j-1}}} 
      -
      e^{C\sum_{h=j-1}^{n-1} \modulo{a_{h+1}-a_h}}
    \right)
    \modulo{\sigma^-}
    \\
    &\leq&
    e^{C\sum_{h=j}^{n-1} \modulo{a_{h+1}-a_h}}
    \left( 
      e^{K_{3}\vert a_j - a_{j-1}\vert}  - e^{C\modulo{a_j - a_{j-1}} }  
    \right)
    \, \modulo{\sigma^-}
    \\
    & \leq &
    (K_{3}-C) \modulo{ a_j - a_{j-1}} \, e^{K_{3}\modulo{a_j - a_{j-1}}}
    \, e^{C\sum_{h=j}^{n-1} \modulo{a_{h+1}-a_h}}
    \, \modulo{\sigma^-} \,.
    \\
    \Delta \Upsilon
    & \leq &
    \left( 
      K_{3} e^{K_{3}\modulo{a^+- a^-}} e^{C\tv(a)} + K_{3} \delta - C
    \right) \,
    \modulo{a_j - a_{j-1}} \, \modulo{\sigma^-}
  \end{eqnarray*}
  and the choice $\delta < 1$ and $C > 2 K_3$ ensures that $\Delta
  \Upsilon < 0$.
\end{proof}

\begin{proofof}{Theorem~\ref{thm:n}}
  First, observe that the construction of the stationary solution
  $\hat u$ directly follows from an iterated application of
  Lemma~\ref{lem:Stationary}. The bound~(\ref{eq:tvhat}) follows from
  the Lipschitz continuity of the map $T$ defined in
  Lemma~\ref{lem:Stationary}.  Define
  \begin{displaymath}
    \tilde{\mathcal{D}}
    =
    \left\{
      u \in
      \hat u +
      \L1 (\reali; A_0)
      \colon u \in \PC \mbox{ and } \Upsilon(u) \leq \delta
    \right\},
  \end{displaymath}
  where $\PC$ denotes the set of piecewise constant functions with
  finitely many jumps. It is immediate to prove that there exists a
  suitable $C_1 > 0$ such that $\frac{1}{C_1} \tv (u)(t) \leq V(t)
  \leq C_1 \tv (u)(t,\cdot)$ for all $(u) \in \tilde{\mathcal{D}}$.
  Any initial data in $\tilde{\mathcal{D}}$ yields an approximate
  solution to~(\ref{eq:HS}) attaining values in $\tilde{\mathcal{D}}$
  by Lemma~\ref{lem:Upsilon}.

  We pass now to the $\L1$-Lipschitz continuous dependence of the
  approximate solutions from the initial datum. Consider two wave
  front tracking approximate solutions $u_1$ and $u_2$ and define the
  functional
  \begin{equation}
    \label{eq:func_cont}
    \Phi\left( u_1,u_2 \right)
    =
    \sum_{j=1}^n \sum_{i=1}^2 \int_0^{+\infty}
    \modulo{s^j_{i}(x)}\, W^j_{i}(x)\, dx\, ,
  \end{equation}
  where $s^j_{i}(x)$ measures the strengths of the $i$-th shock wave
  in the $j$-th pipe at point $x$
  (see~\cite[Chapter~8]{BressanLectureNotes}) and the weights
  $W^j_{i}$ are defined by
  \begin{displaymath}
    W^j_{i}(x)
    =
    1 + \kappa_1 \, A^j_{i}(x) + \kappa_1\, \kappa_2 \,
    \left(
      \Upsilon(u_1)
      +
      \Upsilon(u_2)
    \right)
  \end{displaymath}
  for suitable positive constants $\kappa_1, \kappa_2$ chosen as
  in~\cite[formula~(8.7)]{BressanLectureNotes}.  Here $\Upsilon$ is
  the functional defined in~(\ref{eq:upsilon}), while the $A^j_{i}$
  are defined by
  \begin{eqnarray*}
    A^j_{i}(x)
    & = &
    \sum \left\{ \modulo{\sigma^j_{k_\alpha,\alpha}}
      \colon
      \begin{array}{l}
        x_\alpha<x,\,i<k_\alpha\leq 2
        \\
        x_\alpha>x,\,1\leq k_\alpha<i
      \end{array}
    \right\}
    \\
    & &
    + \left\{
      \begin{array}{ll}
        \displaystyle
        \sum
        \left\{
          \modulo{\sigma^j_{i,\alpha}} \colon
          \begin{array}{l}
            x_\alpha<x,\, \alpha\in \mathcal{J}_j (u_1)
            \\
            x_\alpha>x,\, \alpha\in \mathcal{J}_j (u_2)
          \end{array}
        \right\}
        &
        \textrm{ if }s^j_{i}(x)<0,\vspace{.2cm}
        \\
        \displaystyle
        \sum
        \left\{
          \modulo{\sigma^j_{i,\alpha}} \colon
          \begin{array}{l}
            x_\alpha<x,\, \alpha\in \mathcal{J}_j (u_2)
            \\
            x_\alpha>x,\, \alpha\in \mathcal{J}_j (u_1)
          \end{array}
        \right\}
        &
        \textrm{ if }s^j_{i}(x)\geq0;
      \end{array}
    \right.\nonumber
  \end{eqnarray*}
  see~\cite[Chapter 8]{BressanLectureNotes}. Here, as above,
  $\sigma^j_{i,\alpha}$ is the wave belonging to the $i$-th family,
  sited at $x^\alpha$, with $x^\alpha \in I_j$.  For fixed $\kappa_1$,
  $\kappa_2$ the weights $W^j_{i}(x)$ are uniformly bounded. Hence the
  functional $\Phi$ is equivalent to $\L1$ distance:
  \begin{displaymath}
    \frac1{C_2}\cdot \norma{u_1-u_2}_{\L1}
    \leq
    \Phi\left( u_1,u_2 \right)
    \leq
    \displaystyle
    C_2\cdot \norma{u_1-u_2}_{\L1}
  \end{displaymath}
  for a positive constant $C_2$. The same calculations as
  in~\cite[Chapter 8]{BressanLectureNotes} show that, at any time
  $t>0$ when an interaction happens neither in $u_1$ or in $u_2$,
  \begin{eqnarray*}
    \frac{d}{dt}
    \Phi\left( u_1(t),u_2(t) \right)
    \leq
    C_3\,\epsilon
  \end{eqnarray*}
  where $C_3$ is a suitable positive constant depending only on a
  bound on the total variation of the initial data.

  If $t>0$ is an interaction time for $u_1$ or $u_2$, then, by
  Lemma~\ref{lem:Upsilon}, $\displaystyle \Delta \left[ \Upsilon
    \left(u_1(t)\right) +\Upsilon \left(u_2(t)\right) \right] < 0$
  and, choosing $\kappa_2$ large enough, we obtain
  \begin{displaymath}
    \Delta\Phi\left(u_1(t),u_2(t) \right) < 0\,.
  \end{displaymath}
  Thus, $\Phi \left( u_1(t), u_2(t) \right) - \Phi\left( u_1(s),
    u_2(s) \right) \leq C_2 \, \epsilon \, (t-s)$ for every $0 \leq
  s\leq t$. The proof is now completed using the standard arguments
  in~\cite[Chapter~8]{BressanLectureNotes}.

  The proof that in the limit $\epsilon \to 0$ the semigroup
  trajectory does indeed yield a $\Psi$-solution to~(\ref{eq:HS}) and,
  in particular, that~(\ref{eq:Psii}) is satisfied on the traces, is
  exactly as that of~\cite[Proposition~5.3]{ColomboCorli2}, completing
  the proof of~1.--5.

  Due to the \emph{local} nature of the
  conditions~(\ref{eq:first_int_cond})--(\ref{eq:second_int_cond}) and
  to the finite speed of propagation of~(\ref{eq:HS}), the uniqueness
  conditions~6. and~7.~are proved exactly as in
  Theorem~\ref{thm:Single}.
\end{proofof}

\begin{proofof}{estimate~(\ref{eq:M})}
  We first compute $\partial_{a^+} \Sigma$, with $\Sigma$ defined
  in~(\ref{eq:Sigma}). To this aim, by~2.~in
  Proposition~\ref{prop:Sigma} (in Paragraph 2.3), we may choose
  \begin{displaymath}
    a(x) = \left\{
      \begin{array}{l@{\mbox{ if }}rcl}
        a^- & x & \in & \left] - \infty, -X \right[,
        \\
        \displaystyle
        \frac{a^+ - a^-}{2X} (x+X) + a^- & x & \in & [-X, X]\,,
        \\
        a^+ & x & \in & \left] X , +\infty \right[,
      \end{array}
    \right.
  \end{displaymath}
  so that we may change variable in the integral in~(\ref{eq:Sigma})
  to obtain
  \begin{equation}
    \label{eq:dSigma}
    \partial_{a^+} \Sigma 
    = 
    \partial_{a^+}
    \left( 
      \int_{a^-}^{a^+} p \left( R^a(\alpha, u) \right) \, d\alpha
    \right)
    =
    p(\rho) + O(\Delta a) \,.
  \end{equation}
  Now, estimate~(\ref{eq:M}) directly follows
  inserting~(\ref{eq:pIso}) and~(\ref{eq:dSigma}) in~(\ref{eq:CK})
  and~(\ref{eq:K}).
\end{proofof}

\begin{proofof}{estimates~(\ref{eq:sigma2++})--(\ref{eq:kgrande})}
  Refer to the notation in Figure~\ref{fig:sugiu}, where the pipe's
  section is given by
  \begin{displaymath}
    a(x) 
    =
    \left\{
      \begin{array}{l@{\quad\mbox{ if }}rcl}
        a & x & \in & \left]-\infty, l\right[,
        \\
        a+\Delta a & x &\in &[l, 2l]\,,
        \\
        a & x & \in & \left]2l, +\infty\right[,
      \end{array}
    \right.
  \end{displaymath}
  where $\Delta a>0$. The wave $\sigma_{2}^{+}$ arises from the
  interaction with the first junction and hence
  satisfies~(\ref{eq:sigma2}). Using the pressure law~(\ref{eq:pIso})
  and~(\ref{eq:dSigma}), we obtain
  \begin{eqnarray*}
    \sigma_2^+
    & = &
    \left( 1+\psi\left( u,a\right) \Delta a\right)\sigma_2^- ,
    \qquad \mbox{ where}
    \\
    \psi(a, u)
    & = &
    -\frac{1}{a} \left( 1-\frac{1/2}{1 - (v/c)^{2}}\right) \,.
  \end{eqnarray*}
  Now we iterate the previous bound to estimate the wave
  $\sigma_{2}^{++}$ which arises from the interaction with the second
  junction, i.e.
  \begin{displaymath}
    \sigma_2^{++}
    =
    \left( 
      1 - \psi( a+\Delta a, u^+) \, \Delta a
    \right) \sigma_2^+ ,
  \end{displaymath}
  where, by~(\ref{eq:TT}),
  \begin{displaymath}
    \psi ( a+\Delta a, u^+ )
    =
    \psi\left( 
      a + \Delta a,
      \left( 
        1 + \frac{1}{1-\left( \frac{v}{c}\right) ^{2}} \;\frac{\Delta a}{a}
      \right) 
      \rho,
      \left( 1-\frac{\Delta a}{a}\right) q
    \right).
  \end{displaymath}
  Introduce $\eta = 1/\left(1-(v/c)^{2}\right)$ and $\theta = \Delta a
  /a$ to get the estimate
  \begin{eqnarray*}
    \sigma_2^{++}
    & = &
    \left( 
      1 + \left( 
        \psi( a, u,) -\psi( a+\Delta a, u^+)
      \right)  
      \Delta a
    \right) \sigma_2^-
    \\
    & = &
    \left( 
      1 + 
      \frac{\Delta a}{a} \left( -1+\frac{\eta}{2} \right) 
      +
      \frac{\Delta a}{a+\Delta a} 
      \left( 
        1-\frac{1/2}{1-\left( 
            \frac{1-\theta}{1+\eta\theta}
            \frac{v}{c}
          \right)^{2}}
      \right)
    \right)    
    \sigma_2^-
    \\
    & = &
    \left( 
      1 +
      \frac{\Delta a}{a}
      \left(
        -1 + \frac{\eta}{2} + \frac{1}{1+\theta}
        \left(
          1-\frac{1/2}{1-\left(
              \frac{1-\theta}{1+\theta\eta}
              \frac{v}{c}
            \right)^{2}}
        \right) 
      \right)
    \right)
    \sigma_2^-
  \end{eqnarray*}
  and a further expansion to the leading term in $\Delta a$
  gives~(\ref{eq:sigma2++})--(\ref{eq:kgrande}).
\end{proofof}

\subsection{Proofs Related to Section~\ref{sec:Smooth}}
\label{subs:23}

\begin{proofof}{Lemma~\ref{lem:Weak}}
  If $a \in \C1\left( \reali; [a^-,a^+]\right)$ and $u$ is a weak
  entropy solution of~(\ref{eq:GL}). Then,
  \begin{eqnarray*}
    0
    & = &
    \int_{\reali^+} \int_{\reali} 
    \left( 
      \left[
        \begin{array}{c}
          \rho \\ q
        \end{array}
      \right]
      \partial_t \phi + 
      \left[ 
        \begin{array}{c}
          q
          \\
          P(u)
        \end{array}
      \right]
      \partial_x \phi
      -
      \left[
        \begin{array}{c}
          \frac{q}{a}\partial_{x}a\\ \frac{q^{2}}{a\rho}\partial_{x}a
        \end{array}
      \right]
      \phi
    \right)
    dx \, dt
    \\
    & = &
    \int_{\reali^+} \int_{\reali} 
    \left( 
      \left[
        \begin{array}{c}
          a\rho \\ a q
        \end{array}
      \right]
      \partial_t \frac{\phi}{a}
      + 
      \left[ 
        \begin{array}{c}
          a q
          \\
          a P(u)
        \end{array}
      \right]
      \frac{1}{a}\partial_x \phi
      -
      \left[
        \begin{array}{c}
          a q \\ a \frac{q^{2}}{\rho} 
        \end{array}
      \right]
      \frac{\phi}{a^2} \, \partial_{x}a
    \right)
    dx \, dt
    \\
    & = &
    \int_{\reali^+} \int_{\reali} 
    \left( 
      \left[
        \begin{array}{c}
          a\rho \\ a q
        \end{array}
      \right]
      \partial_t \frac{\phi}{a}
      + 
      \left[ 
        \begin{array}{c}
          a q
          \\
          a P(u)
        \end{array}
      \right]
      \partial_x \frac{\phi}{a}
      +
      \left[
        \begin{array}{c}
          0 \\ p(\rho) \partial_x a
        \end{array}
      \right]
      \frac{\phi}{a}
    \right)
    dx \, dt
  \end{eqnarray*}
  showing that~(\ref{eq:weak}) holds. Concerning the entropy
  inequality, compute preliminarily
  \begin{eqnarray*}
    \nabla(aE(u))
    \left[
      \begin{array}{c}
        \frac{q}{a}\partial_{x}a \\ \frac{q^{2}}{a\rho}\partial_{x}a
      \end{array} 
    \right] 
    & = &
    a \left[ 
      -\frac{q^{2}}{2\rho^{2}}+\int_{\rho}^{\rho_{*}}
      \frac{p(r)}{r^{2}}dr+\frac{p(\rho)}{\rho}
      ,\quad 
      \frac{q}{\rho}
    \right]
    \left[
      \begin{array}{c}
        \frac{q}{a}\partial_{x}a \\ \frac{q^{2}}{a\rho}\partial_{x}a
      \end{array} 
    \right]
    \\
    & = &
    \left( 
      -\frac{q^{3}}{2\rho^{2}}
      +
      q \int_{\rho}^{\rho_{*}}
      \frac{p(r)}{r^{2}}dr+\frac{q}{\rho}p(\rho)
      +
      \frac{q^{3}}{\rho^{2}}\right)\partial_{x}a 
    \\
    & = &
    \frac{q}{\rho} \left( E(u) + p(\rho)\right) \, \partial_{x}a
    \\
    & = &
    F(u) \, \partial_{x}a \,.
  \end{eqnarray*}
  Consider now the entropy condition for~\Ref{eq:GL} and, by the above
  equality,
  \begin{eqnarray*}
    0
    & \leq &
    \int_{\reali^+} \! \int_{\reali}
    \left(
      E(u) \, \partial_t \phi 
      +
      F(u) \, \partial_x \phi -\nabla E(u)
      \left[
        \begin{array}{c}
          \frac{q}{a}\partial_{x}a\\ \frac{q^{2}}{a\rho}\partial_{x}a
        \end{array}
      \right]
      \phi
    \right)
    \, dx \, dt
    \\
    & = &
    \int_{\reali^+} \int_{\reali}
    \Bigl(
    aE(u) \, \partial_t \frac{\phi}{a}+
    aF(u) \, \partial_x \frac{\phi}{a}
    \\
    & & 
    \qquad\qquad
    +
    \left( F(u)\partial_{x}a-\nabla(aE(u))
      \left[
        \begin{array}{c}
          \frac{q}{a}\partial_{x}a\\ \frac{q^{2}}{a\rho}\partial_{x}a
        \end{array}
      \right]\right) \frac{\phi}{a} 
    \Bigl) 
    \, dx \, dt  
    \\
    & = &
    \int_{\reali^+} \int_{\reali}
    \Bigl(
    aE(u) \, \partial_t \frac{\phi}{a}+
    aF(u) \, \partial_x \frac{\phi}{a}
    \\
    & &
    \qquad\qquad+
    \left( F(u)\partial_{x}a-F(u)\partial_{x}a\right) \frac{\phi}{a} 
    \Bigr) 
    \, dx \, dt  
    \\
    & = &
    \int_{\reali^+} \int_{\reali}
    \left(
      aE(u) \, \partial_t \frac{\phi}{a}+
      aF(u) \, \partial_x \frac{\phi}{a}
    \right) 
    \, dx \, dt,
  \end{eqnarray*}
  showing that \Ref{eq:Entropy} holds. The extension to $a \in
  \W{1,1}$ is immediate.
\end{proofof}

\begin{proofof}{Proposition~\ref{prop:Sigma}}
  The regularity condition~\textbf{($\mathbf{\Sigma}$0)} follows from
  the theory of ordinary differential
  equations. Condition~\textbf{($\mathbf{\Sigma}$1)} is immediate.

  Consider now the item~2. If $a_1$ and $a_2$ both
  satisfy~\textbf{(A1)}, are strictly monotone, smooth and have the
  same range, then $a_1 = a_2 \circ \phi$ for a suitable strictly
  monotone $\phi$ with, say $\phi' \geq 0$, the case $\phi' \leq 0$ is
  entirely similar. Note that if $u = \left( R_i(x; u^-) , Q_i (x;
    u^-) \right)$ solves~(\ref{eq:Stationary}) with $a = a_i$, then
  direct computations show that $R_1(x, u^-) = R_2 \left( \phi(x),
    u^-\right)$ and $Q_1(x, u^-) = Q_2 \left( \phi(x),
    u^-\right)$. Hence
  \begin{eqnarray*}
    \Sigma_1 (a^-,a^+; u^-)
    & = &
    \int_{-X}^X p \left( R_1(x; u^-) \right) \, a_1'(x) \, dx
    \\
    & = &
    \int_{-X}^X p \left( R_2\left(\phi(x); u^-\right) \right)
    \, a_2'\left(\phi(x) \right) \, \phi'(x) \, dx
    \\
    & = &
    \int_{-X}^X p \left( R_2(\xi; u^-) \right) \, a_2'(\xi) \, d\xi
    \\
    & = &
    \Sigma_2 (a^-,a^+; u^-) \,.
  \end{eqnarray*}
  Having proved~\textbf{($\mathbf{\Sigma}$0)}
  and~\textbf{($\mathbf{\Sigma}$1)}, we use the map $T$ defined in
  Lemma~\ref{lem:Stationary}. We first prove that $\Sigma$ satisfies
  $\Sigma(a^-,a^+; u^-) + \Sigma \left(a^+, a^- ; T(a^+, a^-; u^-)
  \right) = 0$, given $a$ satisfying~\textbf{(A1)}, strictly monotone
  and with $a(-X) = a^-$, $a(X) = a^+$, let $\tilde a(x) = a^- + a^+ -
  a(x)$. Then, using~2.~proved above, and
  integrating~(\ref{eq:Stationary}) backwards, we have
  \begin{eqnarray*}
    \Sigma \left(a^+, a^-;T(a^-,a^+; u^-) \right)
    & = &
    \int_{-X}^X p\left( \tilde R(x;T(a^-,a^+; u^-)\right) \, \tilde
    a'(x) \, dx
    \\
    & = &
    -  \int_{-X}^X p\left( R(x;a^-,a^+; u^-\right) \, 
    a'(x) \, dx
    \\
    & = &
    -  \Sigma (a^-, a^+; u^-) \,.
  \end{eqnarray*}

  Finally, condition~\textbf{($\mathbf{\Sigma}$2)} follows from the
  the flow property of $R$ and the additivity of the integral. Indeed,
  by~2.~and~3.~we may assume without loss of generality that $a^- <
  a^0 < a^+$. Then, let $q = Q(x; u^-)$ be the $q$ component in the
  solution to~(\ref{eq:Stationary}) with initial condition $u(0) =
  u^-$. Then, if $T$ is the map defined in Lemma~\ref{lem:Stationary},
  we have
  \begin{displaymath}
    T(a^-,a^+; u^-) 
    = 
    \left( 
      R( a^{-1}(a^+);  u^-), Q( a^{-1}(a^+);  u^-)
    \right)
  \end{displaymath}
  so that
  \begin{eqnarray*}
    & &
    \Sigma(a^-,a^+; u^-)
    \\
    & = &
    \int_{-X}^{X} p \left( R(x,  u^-) \right) \, a'(x) \, dx
    \\
    & = &
    \int_{-X}^{a^{-1}(a^0)} p \left( R(x,  u^-) \right) \, a'(x)
    \, dx
    +
    \int_{a^{-1}(a^0)}^{X} p \left( R(x,  u^-) \right) \, a'(x) \, dx
    \\
    & = &
    \int_{-X}^{a^{-1}(a^0)} p \left( R(x,  u^-) \right) \, a'(x) \, dx
    \\
    & &
    +
    \int_{a^{-1}(a^0)}^{X} p 
    \left( 
      R\left(
        x, R(a^{-1}(a^0), u^-), Q(a^{-1}(a^0), u^-) 
      \right)
    \right) \, a'(x) \, dx
    \\
    & = &
    \Sigma(a^-,a^0; u^-) 
    +
    \Sigma \left( a^0,a^+;T(a^-,a^+; u^-) \right)
  \end{eqnarray*}
  proving~1.
\end{proofof}

\begin{proofof}{Theorem~\ref{thm:W11}}
  Fix $\bar a>0$, and $\bar u \in A_0$. Choose $M,\Delta ,L,\delta$ as
  in Theorem~\ref{thm:n}. With reference to these quantities, let $a$
  satisfy~\textbf{(A1)}. For $n \in \naturali$, let $a_n, \alpha_n,
  \beta_n$ be as in Proposition~\ref{prop:an}.
  \begin{figure}[htpb]
    \centering
    \begin{psfrags}
      \psfrag{mX}{$-X$} \psfrag{X}{$X$} \psfrag{x1}{$x^j_n$}
      \psfrag{y1}{$y^j_n$} \psfrag{y2}{$y^{j+1}_n$} \psfrag{x}{$x$}
      \includegraphics[width=9cm]{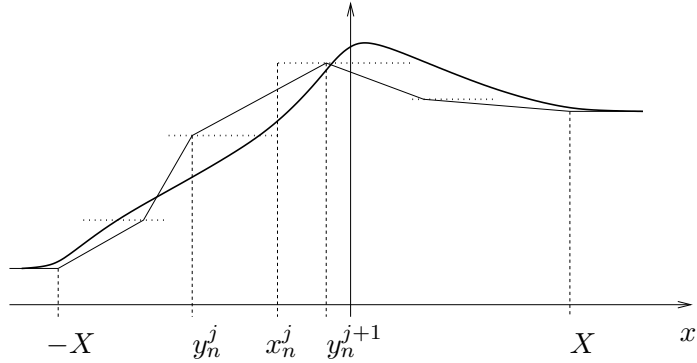}
    \end{psfrags}
    \caption{The thick line is the graph of $a=a(x)$, the dotted line
      represents $a_n$ while the polygonal line is $\alpha_n$}
    \label{fig:alpha}
  \end{figure}
  Note that $\alpha_n$ is piecewise linear and
  continuous. 
  By~(iii), we have that $\alpha_n \to a$ and $a_n \to a$ in
  $\L1$. Moreover, $\tv(\alpha_n) \leq M$ and $\tv(a_n) \leq M$ and,
  for $n$ sufficiently large, $a_n(\reali) \subseteq \left]\bar a -
    \Delta, \bar a + \Delta \right[$. Hence, for $n$ large, $a_n$
  satisfies~\textbf{(A0)}. Call $S^n$ the semigroup constructed in
  Theorem~\ref{thm:n} and denote by $\mathcal{D}^n$ its domain.

  Let $u_{n}^{0}$ be a sequence of initial data in
  $\mathcal{D}^n$. The $S^n$ are uniformly Lipschitz in time and
  $S^n_tu_{n}^{0}$ have total variation in $x$ uniformly bounded in
  $t$. Hence, by~\cite[Theorem~2.4]{BressanLectureNotes}, a
  subsequence of $u_{n}(t) = S^n_tu_{n}^{0}$ converges pointwise
  a.e.~to a limit, say, $u$. For any $\phi \in \Cc1
  (\pint\reali^+\times \reali; \reali)$ and for any fixed $n$, let
  $\epsilon >0$ be sufficiently small and introduce a
  $\Cc\infty(\reali;\reali)$ function $\eta_\epsilon$ such that
  \begin{displaymath}
    \begin{array}{c@{\quad\mbox{ for all }\quad}rcl}
      \eta_\epsilon(x) = 0 
      & 
      x & \in & 
      \bigcup_{j=1}^{m_n-1} [x^j_n-\epsilon, x^j_n+\epsilon] \,,
      \\[5pt]
      \eta_\epsilon(x) = 1
      & 
      x & \in &
      \bigcup_{j=1}^{m_n-2} [x^j_n+2\epsilon, x^{j+1}_n-2\epsilon] \,.
    \end{array}
  \end{displaymath}
  Thus, we have
  \begin{eqnarray*}
    & &
    \int_{\reali^+} \! \int_{\reali}  \!
    \left( 
      \left[ \!
        \begin{array}{c}
          a_n\rho_n \\ a_n q_n
        \end{array}
        \! \right]
      \partial_t \phi + 
      \left[ \!
        \begin{array}{c}
          a_n q_n
          \\
          a_n P(u_{n})
        \end{array}
        \! \right]
      \partial_x \phi
    \right)\!
    {dx} \, {dt}
    \\
    & = &
    \lim_{\epsilon\to 0}
    \int_{\reali^+} \! \int_{\reali}  \!
    \left( 
      \left[ \!
        \begin{array}{c}
          a_n\rho_n \\ a_n q_n
        \end{array}
        \! \right]
      \eta_\epsilon \,
      \partial_t \phi + 
      \left[ \!
        \begin{array}{c}
          a_n q_n
          \\
          a_n P(u_{n})
        \end{array}
        \! \right]
      \eta_\epsilon \,
      \partial_x \phi
    \right)\!
    {dx} \, {dt}
    \\
    & = &
    \lim_{\epsilon\to 0}
    \int_{\reali^+} \! \int_{\reali}  \!
    \left( 
      \left[ \!
        \begin{array}{c}
          a_n\rho_n \\ a_n q_n
        \end{array}
        \! \right]
      \partial_t (\eta_\epsilon \,\phi )+ 
      \left[ \!
        \begin{array}{c}
          a_n q_n
          \\
          a_n P(u_{n})
        \end{array}
        \! \right]
      \partial_x (\eta_\epsilon \,\phi)
    \right)\!
    {dx} \, {dt}
    \\
    & &
    \qquad
    -
    \lim_{\epsilon\to 0}
    \int_{\reali^+} \! \int_{\reali}  \!
    \left[ \!
      \begin{array}{c}
        a_n q_n
        \\
        a_n P(u_{n})
      \end{array}
      \! \right]
    \phi \, \partial_x \eta_\epsilon \,
    {dx} \, {dt} \,.
  \end{eqnarray*}
  The first summand in the latter term above vanishes by
  Definition~\ref{def:WeakPsi} applied in a neighborhood of each
  $x^j_n$. The second summand, by the $\BV$ regularity of $u_n$,
  converges as follows:
  \begin{eqnarray*}
    & &
    -
    \int_{\reali^+} \! \int_{\reali}  \!
    \left( 
      \left[ \!
        \begin{array}{c}
          a_n\rho_n \\ a_n q_n
        \end{array}
        \! \right]
      \partial_t \phi + 
      \left[ \!
        \begin{array}{c}
          a_n q_n
          \\
          a_n P(u_{n})
        \end{array}
        \! \right]
      \partial_x \phi
    \right)\!
    {dx} \, {dt}
    \\
    & = &
    \lim_{\epsilon\to 0}
    \int_{\reali^+} \! \int_{\reali}  \!
    \left[ \!
      \begin{array}{c}
        a_n q_n
        \\
        a_n P(u_{n})
      \end{array}
      \! \right]
    \phi \, \partial_x \eta_\epsilon \,
    {dx} \, {dt}
    \\
    & = &
    \sum_{j=1}^{m_n-1}
    \int_{\reali^+} \left[
      \begin{array}{c}
        a_n(x^j_n+) q_n(x^j_n+) - a_n(x^j_n-) q_n(x^j_n-)
        \\
        a_n(x^j_n+) P_n(x^j_n+) - a_n(x^j_n-) P_n(x^j_n-)
      \end{array}
    \right]
    \phi(t, x^j_n) \, {dt}
    \\
    & = &
    \sum_{j=1}^{m_n-1}
    \int_{\reali^+} \left[
      \begin{array}{c}
        0
        \\
        \Sigma\left( 
          a_n(x^j_n-), a_n(x^j_n+),u(t,x^j_n-) 
        \right)
      \end{array}
    \right]
    \phi(t, x^j_n) \, {dt} \,.
  \end{eqnarray*}
  We proceed now considering only the second component. Using the map
  \begin{eqnarray*}
    \phi_n (t,x) 
    & = &
    \phi(t,x) \, \chi_{\strut ]-\infty, y^1_n[}(x) 
    + 
    \sum_{j=1}^{m_n-1} \phi(t, x^j_n) \, \chi_{[y^j_n, y^{j+1}_n [}(x)
    \\
    & &
    + \phi(t,x) \, \chi_{\strut ]y^{m_n}_n,  +\infty[}(x) \,,
  \end{eqnarray*}
  we obtain
  \begin{eqnarray*}
    \! & & \!\!\!
    \sum_{j=1}^{m_n-1}
    \int_{\reali^+} 
    \Sigma\left( 
      a_n(x^j_n-), a_n(x^j_n+),u(t,x^j_n-) 
    \right)
    \phi(t, x^j_n) \, {dt}
    \\
    \! & = & \!\!\!
    \sum_{j=1}^{m_n-1}
    \int_{\reali^+} 
    \Sigma\left( 
      a_n(y^j_n), a_n(y^{j+1}_n),u(t,x^j_n-) 
    \right)
    \phi(t, x^j_n) \, {dt}
    \\
    \! & = & \!\!\!
    \sum_{j=1}^{m_n-1}
    \int_{\reali^+} 
    \Sigma\left( 
      \alpha_n(y^j_n), \alpha_n(y^{j+1}_n),u(t,x^j_n-) 
    \right)
    \phi(t, x^j_n) \, {dt}
    \\
    \! & = & \!\!\!
    \sum_{j=1}^{m_n-1}
    \int_{\reali^+} 
    \int_{y^j_n}^{y^{j+1}_n} \!\!\!\!
    p \! \left(
      R^{\alpha_n} \left(x;u_n(t,x^j_n-)\right)
    \right)
    \alpha'_n(x) \, dx \; 
    \phi(t, x^j_n) \, {dt}
    \\
    \! & = & \!\!\!
    \int_{\reali^+} 
    \sum_{j=1}^{m_n-1}
    \int_{y^j_n}^{y^{j+1}_n} \!\!\!\!
    p \! \left(
      R^{\alpha_n} \left(x;u_n(t,x^j_n-)\right)
    \right)
    \alpha'_n(x) \, dx \; 
    \phi(t, x^j_n) \, {dt}
    \\
    \! & = & \!\!\!
    \int_{\reali^+} 
    \int_{\reali}
    \sum_{j=1}^{m_n-1}
    p \! \left(
      R^{\alpha_n} \left(x;u_n(t,x^j_n-)\right)
    \right)
    \\
    & &
    \qquad\qquad\qquad\qquad\qquad\qquad
    \times
    \alpha'_n(x) \, 
    \phi(t, x^j_n) \, 
    \chi_{[y^j_n, y^{j+1}_n [}(x) \,
    dx \, {dt}
    \\
    \! & \to & \!\!\!
    \int_{\reali^+} \int_{\reali} 
    p\left( \rho(x) \right) \,
    \partial_x a(x) \,
    \phi(t, x^j_n) \, {dx} \, {dt}
    \qquad \mbox{ as } n \to +\infty \,,
  \end{eqnarray*}
  where we used~(i) in the choice of the approximation $\alpha_n$.

  We thus constructed a solution to~(\ref{eq:GL2}), for any initial
  datum in $\mathcal{D}$. Note that this solution
  satisfies~(\ref{eq:first_int_condBIS})--(\ref{eq:second_int_condBIS}),
  as can be proved using exactly the techniques
  in~\cite[Theorem~8]{GuerraMarcelliniSchleper}. Therefore, the whole
  sequence $u_n$ converges to a unique limit $u$, which is Lipschitz
  with respect to time. This uniqueness implies the semigroup
  property~2.~in Theorem~\ref{thm:W11}. The Lipschitz continuity with
  respect to the initial datum follows from the uniform Lipschitz
  regularity of the approximate solutions $u_n$, completing the proof
  of 3. Finally, 6.~is proved exactly as
  in~\cite[Theorem~8]{GuerraMarcelliniSchleper}.
\end{proofof}

{\small{

    \bibliographystyle{abbrv}

    \bibliography{ColomboMarcellini_revised}

  }}

\end{document}